%% file: main.tex
\documentclass[preprint,3p]{elsarticle}  %
\biboptions{sort}
\input{preambles/general_preamble}

\input{preambles/crdme_preamble}

\renewcommand{\vec}[1]{\mathbf{#1}}

\newcommand{\densityletter}{f}
\newcommand{\vecdensity}{\vec{\densityletter}}
\newcommand{\density}[1]{\densityletter^{#1}}

\newcommand{\probletter}{\MakeUppercase{\densityletter}}
\newcommand{\veccdf}{\vec{\probletter}}
\newcommand{\cdfvec}{\veccdf} 
\newcommand{\cdf}[1]{\probletter^{#1}}

\newcommand{\cdfeq}{\bar{\probletter}^{\specvec}}

\newcommand{\speciesletter}{s}
\newcommand{\specvec}{\vec{\speciesletter}}
\newcommand{\scounts}[1]{\speciesletter_{#1}}
\newcommand{\plocsletter}{q}
\newcommand{\plocs}{\vec{\plocsletter}}
\newcommand{\nspecies}{M}
\newcommand{\slocs}{(\plocs^{\scounts{1}},\ldots,\plocs^{\scounts{\nspecies}})}
  
\newcommand{\denseq}{\bar{\densityletter}^{\specvec}}
\newcommand{\denseqplus}{\bar{\densityletter}^{\specvec^+}}

\newcommand{\Rabc}{R^{(a,b,c)}}
\newcommand{\slocsabc}{(\plocs^{a},\plocs^{b},\plocs^{c})}
\newcommand{\fwdop}[3]{\Lambda^+_{#3 | #1,#2}}
\newcommand{\annihilation}[1]{\mathcal{A}_{#1}}

\newcommand{\bwdop}[3]{\Lambda^-_{#1,#2|#3}}
\newcommand{\hopop}[3]{\Lambda^{#1}_{#3|#2}}
\newcommand{\node}[1]{\mathbf{x}_{#1}}

\newcommand{\transport}[1]{L^{\mathbf{#1}}}
\newcommand{\reaction}[1]{R^{\mathbf{#1}}}
\newcommand{\boldq}{\mathbf{q}}
\newcommand{\bolds}{\mathbf{s}}
\newcommand{\bq}{\mathbf{q}}
\newcommand{\bs}{\mathbf{s}}
\newcommand{\bd}{\mathbf{d}}
\newcommand{\bx}{\mathbf{x}}

\newcommand{\partialt}[1]{\frac{\partial #1}{\partial t}}
\newcommand{\derivt}[1]{\frac{d #1}{d t}}

\begin{document}

\begin{frontmatter}

\title{A Convergent Reaction-Drift-Diffusion Master Equation with Interaction Potentials}

\author[bu]{Max Heldman}
\ead{heldmanm@bu.edu}
\author[bu]{Samuel A. Isaacson}
\address[bu]{Department of Mathematics and Statistics, Boston University, Boston, MA}
\ead{isaacsas@bu.edu}

\begin{abstract}
We extend the convergent reaction-diffusion master equation (CRDME) framework on unstructured grids to include two-body interaction potentials between particles. Such interactions are essential for modeling systems in which particles experience electrostatic repulsion, volume exclusion, or van der Waals attraction. Our approach discretizes the weak form of the forward equation for a particle-based stochastic reaction-drift-diffusion model to obtain a new CRDME by combining two complementary techniques: an edge-averaged finite element method with mass-lumping quadrature for spatial transport terms that include potential interactions, and a mix of finite volume and mass-lumping quadrature approximations for reaction terms. This yields transition rates for a spatial jump process of particles hopping between dual mesh voxels and undergoing reactions. The resulting spatially discrete CRDME retains key physical properties of the continuum system in closed domains. It preserves the Gibbs-Boltzmann equilibrium distribution in the absence of reactions, and satisfies detailed balance of both drift-diffusion and reactive fluxes for systems with reversible reactions. Numerical comparisons with Brownian dynamics simulations validate our method, and we demonstrate empirical second-order convergence as the mesh is refined.
\end{abstract}


\end{frontmatter}

\newtheorem{remark}{Remark}
\graphicspath{{manuscript/img/}}
\input{manuscript/introduction}

\input{manuscript/background}

\input{manuscript/results}

\input{manuscript/derivation-general}

\section{Conclusion}\label{section:conclusion}
We have extended the convergent reaction-diffusion master equation (CRDME) framework to incorporate two-body interaction potentials between particles, significantly broadening the class of systems that can be accurately simulated. Our derived hopping and reaction rates ensure that the discrete system maintains detailed balance for reversible reactions while preserving the Gibbs-Boltzmann equilibrium distribution. The numerical experiments demonstrate empirical second-order convergence rates as the mesh spacing is reduced, and show good agreement with Brownian dynamics simulations. Our work provides a systematic numerical framework for spatially discrete stochastic simulations of interacting particle systems, combining numerical accuracy with thermodynamic consistency. 

Past CRDME literature (e.g., \cite{isaacson_convergent_2013,isaacson_unstructured_2018}) includes empirical evidence that statistics of the CRDME for the reaction-diffusion problem (using \eqref{eq:discrete_fwd_standard} and \eqref{eq:discrete_bwd_standard} for the reaction rates) converge to their continuum analogues at $O(h^2)$. More recently, \cite{isaacson_zhang_crddme_2024} does the same for the case of one-body potentials discretized using EAFE, using a detailed-balance-satisfying modification of \eqref{eq:discrete_fwd_standard} and \eqref{eq:discrete_bwd_standard} for reactions. Our results in Section~\ref{section:numerics_pair_potentials} suggest that our new CRDME with interaction potentials maintains this rate of convergence.

Convergence of the density governed by \eqref{eq:discrete-FKE-full} to that governed by \eqref{eq:FKE-full}, or even of CRDME statistics to their continuum analogues, has not been proven rigorously. However, various works have considered the accuracy of $P^1$ finite element discretizations of related PDEs using one or more of the approximations mentioned above. In what follows, by a reaction-diffusion or reaction-drift-diffusion PDE we mean a linear PDE containing zeroth-order (nonderivative) terms analogous to the reaction terms in the forward equation~\eqref{eq:FKE-full}. We highlight a few main results here:
\begin{itemize}
\item For scalar, time-dependent reaction-diffusion PDEs, applying mass-lumping quadrature to the time derivative term is standard, and, e.g., \cite{thomee_galerkin_2006} shows $O(h^2)$ error estimates in the $L^\infty([0,T]; L^2(\Omega))$-norm for such problems.
\item For scalar, time-independent reaction-drift-diffusion PDEs with piecewise $W^{2,\infty}(\Omega)$ reaction coefficients, \cite{Heldman_2024} shows $O(h^2)$ error estimates in the $L^2(\Omega)$-norm when \eqref{eq:discrete_fwd_standard} and \eqref{eq:discrete_bwd_standard} are applied to the reaction term ---  see also the doctoral thesis \cite{Heldman_2023}, Chapter 5, for additional analytical work on the accuracy of \eqref{eq:discrete_fwd_standard}--\eqref{eq:discrete_bwd_smooth} for both smooth and discontinuous reaction kernels.
\item For scalar, time-independent reaction-drift-diffusion PDEs with sufficiently smooth coefficients, \cite{xu_monotone_1999} shows $O(h)$ estimates in the $H^1(\Omega)$-norm when the EAFE discretization is applied to the transport operator and mass-lumping quadrature to the zeroth-order term --- additional empirical (e.g., \cite{miller_scharfetter-gummel_1991}) and analytical \cite{Heldman_2024} work suggests that the resulting discretization is second-order accurate in the $L^2(\Omega)$-norm, at least in some cases.
\end{itemize}
We expect that carefully applying elements of the various estimates listed above to our discretization of \eqref{eq:FKE-full} with the techniques of \cite{raviart_use_1975,thomee_galerkin_2006} can lead to convergence rate estimates for the discretized solution. We plan to address the question of convergence in detail analytically in future work. 

Beyond the question of the rigorous rate of convergence in the mesh width, a number of further extensions of the current work are possible. For example, this work has focused on deriving an accurate and convergent method; in future work we hope to report on optimized exact SSAs that can be used to efficiently simulate the resulting CRDME with interaction potentials. We have also found that the use of rejection methods can significantly simplify the implementation of CRDME SSAs, for example bypassing the need to precalculate reactive transition rates, while for many parameter regimes keeping overall simulation time close to implementations that precalculate rates~\cite{Heldman_2023}.

\section*{Acknowledgments}
The work of both authors was supported by the Army Research Office under grant W911NF2010244 and by the National Science Foundation under grant NSF-DMS 1902854. The work of SAI was additionally supported by Army Research Office grant W911NF2510078 and National Science Foundation  grant NSF-DMS 2325185.

\section*{Declaration of generative AI and AI-assisted technologies in the writing process}
Beginning with the editing of a later, comprehensive draft of this manuscript, SAI used GitHub Copilot and Claude Code as part of the Visual Studio Code editor. This included auto-completion suggestions for text, and suggestions for revisions to written text to improve clarity of exposition and grammar.

SAI also subsequently used GPT 5.6 Sol on Max effort to review a near-final draft of the manuscript for clarity, wording, and correctness. It suggested several edits to fix typos and make minor corrections to equations, suggested focusing the discussion at the end of Section~\ref{section:derivations_reactions_multiparticle} on what happens within stoichiometric compatibility classes, and suggested the argument in the second-to-last paragraph of Section~\ref{section:derivations} (i.e., that, by constructing an equilibrium on each stoichiometric compatibility class and invoking finiteness and irreducibility of the class-restricted CRDME, one can conclude that for any given initial distribution one converges to a uniquely determined equilibrium distribution satisfying detailed balance).

MH used Claude Code as part of the Visual Studio Code editor for final edits, including grammar, consistency of notation and grammatical conventions, and clarity of the writing. Claude Code also identified several minor errors in mathematical formulae that needed to be corrected. After using these tools, the authors reviewed and edited the content as needed and take full responsibility for the final content of the published article.

\bibliographystyle{elsarticle-num}
\bibliography{manuscript/references}

\end{document}

%% file: preambles/general_preamble.tex
\usepackage{amssymb}
\usepackage{amsmath}
\usepackage{mathtools}
\usepackage{graphicx}
\usepackage{tikz}
\usepackage{pgfplots}
\usepackage{url}
\usepackage[colorlinks=true,linkcolor=blue,citecolor=blue,urlcolor=blue]{hyperref}

\newcommand{\NN}{\mathbb{N}}

\newcommand{\RR}{\mathbb{R}}

\newcommand{\emptytext}[1]{ \mathbin{\color{white}#1}}

\renewcommand{\vec}[1]{\boldsymbol{#1}}

\newcommand{\paren}[1]{\left(#1\right)}
\newcommand{\brac}[1]{\left[#1\right]}
\newcommand{\abs}[1]{\left|#1\right|}

\DeclareMathOperator{\Span}{\textrm{span}}

%% file: preambles/crdme_preamble.tex
\newcommand{\qvec}{{\mathbf{q}}}
\newcommand{\qplus}{\mathbf{q}^{+}}
\newcommand{\qminus}{\mathbf{q}^{-}}

\newcommand{\multinode}[2]{\mathbf{#1}^{\mathbf{#2}}}

\newcommand{\E}{\mathbb{E}}
\newcommand{\avg}[1]{\E\brac{#1}}

\newcommand{\bi}{\vec{i}}
\newcommand{\bj}{\vec{j}}
\newcommand{\bk}{\vec{k}}

\newcommand{\vq}{\vec{q}}

\newcommand{\vqa}{\vec{q}^{a}}
\newcommand{\vqb}{\vec{q}^{b}}
\newcommand{\vqc}{\vec{q}^{c}}

\def\R{\mathbb{R}}

\newcommand{\Kd}{K_{\text{eq}}}

\renewcommand{\epsilon}{\varepsilon}

\newcommand{\qa}{\mathbf{q}^a}
\newcommand{\qb}{\mathbf{q}^b}
\newcommand{\qc}{\mathbf{q}^c}

%% file: manuscript/introduction.tex
\section{Introduction}\label{section:introduction}

Particle-based stochastic reaction-diffusion (PBSRD) systems model individual molecules as they react and diffuse randomly in a physical domain $\Omega\subseteq \RR^d$, where $d\in \{1,2,3\}$. PBSRD systems track the random motion and interactions of individual molecules instead of their concentrations, as in macroscopic reaction-diffusion partial differential equations (PDEs), but do not consider more refined models of molecular motion which propagate trajectories of individual atoms via Newton's laws of motion, as in microscopic molecular dynamics (MD) simulations \cite{shaw_millisecond-scale_2009}. PBSRD systems are therefore considered mesoscopic models of particle dynamics, typically incorporating little (explicit) modeling of the internal structure of the molecules while maintaining their individuality.

Despite being more macroscopic than MD models, PBSRD simulators can include physics at varying levels of detail. PBSRD simulators that can handle complex geometries and include potentials between pairs of particles, for example modeling soft-core repulsion or electrostatic interactions, are considered among the most detailed. Such two-body interactions have been used to model a variety of biophysical phenomena, including molecular crowding in T cell signaling~\cite{Siokisetal2018}, protein aggregation~\cite{WinkelmannMckeanVlasov2026,Carrillo2019}, and electrostatic interactions between charged biomolecules~\cite{Huber2019}. Our goal in this work is to develop a convergent jump process approximation of PBSRD models that include two-body interaction potentials via spatial discretization of the underlying forward equation. We would like such an approximation to enable the efficient numerical solution of such systems, while also preserving key physical properties of the underlying continuum model.

Our approach will be based on extending the convergent reaction-diffusion master equation (CRDME) framework developed in \cite{isaacson_convergent_2013,isaacson_unstructured_2018,zhang_particle-based_nodate} to include two-body interaction potentials. Originally conceived in \cite{isaacson_convergent_2013} as a lattice-based, confined-diffusion discretization of spatially continuous PBSRD models on structured meshes, the CRDME has since been extended to complex geometries \cite{isaacson_unstructured_2018}, and to include additional physics such as one-body potentials \cite{isaacson_zhang_crddme_2024}. We show that our new CRDME with interaction potentials retains a number of the advantages cited in the research above. These include exhibiting (empirical) second-order convergence rates to the continuum solution as the mesh is refined, and exact preservation of physical properties, such as detailed balance for reversible reactions and the Gibbs-Boltzmann equilibrium distribution.

To derive our numerical method, we take the perspective of, e.g., \cite{isaacson_convergent_2013,isaacson_unstructured_2018,isaacson_mean_2022,isaacson_zhang_crddme_2024}, and treat the PBSRD system as a continuum model which can be represented as a measure-valued stochastic process (MVSP). This MVSP encodes the dynamics of the state of the system, equivalently representing it as singular concentration fields or as the number of particles of each species along with their positions (the representation we will subsequently use). The probability density that the MVSP has a given state can be represented by its forward Kolmogorov equation, a potentially arbitrarily high-dimensional system of partial integro-differential equations (PIDEs) in Fock space. An appropriate discretization of this forward Kolmogorov equation corresponds to a master equation (the CRDME), i.e., the forward equation for a continuous-time Markov chain (CTMC) representing a spatial jump process. While the CRDME itself is generally too high dimensional to solve using ordinary differential equation (ODE) integrators, the associated CTMC can be exactly sampled using any variant of Gillespie's stochastic simulation algorithm (SSA) \cite{gillespie_general_1976} (also known as kinetic Monte Carlo or Doob's method). The CRDME discretization can also be tailored to satisfy the \textit{mathematical} detailed balance relation for CTMCs, which coincides with \textit{physical} detailed balance conditions for the discretized system. 

The CRDME is a purely \textit{spatial} discretization of the forward equation, dividing $\Omega$ into a lattice of polygons or polyhedra, commonly called voxels, and treating particles in individual voxels as indistinguishable. The CTMC is simulated as a jump process where particles may hop between voxels, simulating diffusive transfer, and react with particles in nearby voxels. At this coarse level of description, the CRDME is similar to the reaction-diffusion master equation (RDME) from which it inherits its name. The difference between the two methods lies in how the transition rates are derived: while the two methods typically treat diffusion in the same way, the RDME only allows particles in the same voxel to react, whereas the CRDME allows reactions between particles within a \textit{geometric} radius (for a Doi interaction). This means that the RDME only gives an accurate representation of the reaction physics for certain mesh widths, which may or may not be sufficiently small to also resolve other physical processes in the model (such as diffusion or other reactions) \cite{isaacson_reaction-diffusion_2009,IsaacsonRDMELimsII,Hellander:2012jk}. The CRDME is designed with the goal that simulation statistics converge to their continuum analogues in the underlying PBSRD model as the mesh size tends to zero, i.e., there is a mesh size which can accurately resolve all physical processes from the continuum model. 

An alternative simulation approach is given by the Brownian dynamics method, e.g., as implemented in the ReaDDy software \cite{schoneberg_readdy_2013,hoffmann_readdy_2019}, which instead discretizes the underlying PBSRD continuum model in time to obtain an approximating spatially continuous discrete-time stochastic process. Statistics computed for the same physical model by our CRDME and the Brownian dynamics method agree as the mesh spacing in the CRDME and the Brownian dynamics time step are reduced, as we confirm in Section~\ref{section:numerics_pair_potentials}. The ReaDDy software allows reactions of order $\leq 2$ and particle transport by diffusion with one- and two-body potentials. Moreover, in~\cite{frohner_reversible_2018} the authors of ReaDDy introduce modified reaction rates that ensure detailed balance for systems of reversible reactions with interaction potentials, and develop a practical acceptance-rejection scheme for such reactions which approximately (up to temporal discretization error of the underlying stochastic differential equation solver) preserves detailed balance. We borrow from \cite{frohner_reversible_2018} in formulating our detailed-balance-satisfying PBSRD continuum model.

The CRDME discrete-space, continuous-time discretization has several properties that distinguish it from the Brownian dynamics discrete-time, continuous-space discretization. The finite element method (FEM) spatial discretization that we use in this work allows one to naturally handle complex geometries \cite{isaacson_unstructured_2018} and even manifolds, and provides the scaffolding to conduct rigorous error analysis. The structured representation as a CTMC arising via spatial discretization of the PBSRD model allows us to exactly preserve detailed balance of both drift-diffusion and reactive fluxes in the CRDME, and to exactly preserve the Gibbs-Boltzmann equilibrium distribution.

The main difficulties we encounter in deriving the new CRDME are due to the particle motion and reactions now depending on the location of all other particles in the system. Therefore, we need to directly discretize the full multiparticle PBSRD model \eqref{eq:FKE-full}, and obtain hopping and reactive transition rates that are functions of the positions of all particles. This is in contrast to earlier CRDME studies~\cite{isaacson_convergent_2013,isaacson_unstructured_2018,isaacson_zhang_crddme_2024}, where it was sufficient to discretize simpler one- and two-particle models to obtain transition rates that end up being the same as what one obtains by directly discretizing the general multiparticle system. We discretize the general PBSRD model by using a tensor product finite element basis on particle configuration space, combined with the edge-averaged finite element (EAFE) method and mass-lumping quadrature for the spatial transport operator, and a mix of mass-lumping and finite-volume-type quadratures for the reaction operator. For many applications, especially in cell biology, interaction potentials are short range, so that one usually only has to resolve interactions between particles that are within several mesh voxels of each other in simulations of the jump process associated with the CRDME. Our discretization approach results in spatial hopping transition rates that are efficient to compute and use in SSA simulations for systems involving short-range interactions. 

The rest of the paper is organized as follows. In the remainder of this section, we introduce the general form of the continuum model and a model problem with interaction potentials and a reversible reaction. In Section~\ref{section:simulation_algorithm_background}, we introduce concepts related to our spatial discretization, discussing previous work on deriving CRDME reactive transition rates and transport rates in systems with no potentials and one-body potentials, respectively. In Section~\ref{section:simulation_algorithm_pair_potentials}, we present our main results, the CRDME hopping and reaction rates for reaction-diffusion systems with interaction potentials. In Section~\ref{section:numerics_pair_potentials}, we give numerical results for the new method, both demonstrating empirical convergence rates and validating against the ReaDDy Brownian dynamics software. In Section~\ref{section:derivations}, we give detailed derivations of the transition rates presented in Section~\ref{section:simulation_algorithm_pair_potentials}. For the general multiparticle reversible $A + B \rightleftharpoons C$ reaction in a closed system, we also establish consistency of the CRDME with Gibbs-Boltzmann spatial configurations at equilibrium, and establish that it supports detailed balance of the resulting drift-diffusion and reactive fluxes. Finally, we conclude by discussing some potential future directions.

\subsection{Notation and setup}\label{section:introduction_model_problem}

We begin with our continuum model, in which individual molecules of types $\mathbf{S} = (S_1,S_2,\dots,S_M)$ and copy numbers $\mathbf{s} = (s_1,s_2,\dots,s_M) \in \NN^M_0$ are represented as points in a closed, bounded, and connected domain $\Omega\subseteq \RR^d$. We treat $\bs$ as a multi-index, with the standard definitions for operations on it, e.g., $\bs! = \prod_{i = 1}^M s_i!$. The particles can experience drift, diffuse, and react randomly according to rates depending on the chemical species involved and the current state of the system. We then discretize the continuum system by subdividing $\Omega$ into a set of polygonal voxels $\{V_i\}_{i=1}^N$ and allowing molecules to move by an approximating spatial jump process, hopping between voxels and reacting with derived transition rates. We determine the hopping and reactive transition rates by a convergent discretization of the forward Kolmogorov equation, written using Doi's notation \cite{doi_second_1976}  as
\begin{equation}\label{eq:FKE-full}
\begin{split}
\partialt{\density{\specvec}}(\plocs,t) &= [\transport{\specvec} \density{\specvec}] (\plocs,t) + [\reaction{\specvec}  \vecdensity](\plocs,t),\\
\density{\specvec}(\boldq,0) &= \density{\specvec}_0(\plocs),
\end{split}
\end{equation} 
of the continuum particle system for the density $\vecdensity = [\density{\specvec}]_{\specvec\in (\NN_0)^\nspecies}$, with components $\density{\specvec}:\Omega^{|\specvec|}\times[0,T]\to \RR$ for some final time $T > 0$. The quantity 
\begin{equation}\label{eq:particle_count_probabilities}
\pi_{\specvec}(t) = \frac{1}{\specvec !}\int_{\Omega^{|\specvec|}} \density{\specvec}(\plocs,t)d\plocs  
\end{equation} 
represents the probability of observing a given number $\bolds$ of each particle type at time $t$, so that after normalization $\frac{1}{\bs! \pi_{\specvec}(t)}\density{\specvec}(\boldq,t)$ gives the probability density for the spatial particle configuration 
\begin{align*}\plocs = \slocs \in \Omega^{|\specvec|} = \Omega^{s_1} \times \Omega^{s_2} \times \dots \times \Omega^{s_M}\text{ for }\plocs^{s_m} := (q^{s_m}_1,\dots,q^{s_m}_{s_m}), 
\end{align*} 
conditional on the particle count $\mathbf{s}$. Here $\plocs^{s_m}$ represents the positions of particles of species $S_m$, with $\plocs$ the full state vector for the position of all particles. We always assume a no-flux boundary condition for \eqref{eq:FKE-full} unless stated otherwise.

The evolution of the density is governed by a drift-diffusion transport operator $L^{\mathbf{s}}$ and a reaction operator $R^{\mathbf{s}}$. In this paper, we assume the particular form 
\begin{equation*}
L^{\mathbf{s}}\density{\specvec}(\mathbf{q},t) = \nabla_{\mathbf{q}} \cdot \left[
\paren{\begin{smallmatrix}
    D^{S_1} I^{s_1d} & &  & \\
    &       & \ddots & \\
    & & &  D^{S_M} I^{s_Md}
\end{smallmatrix}} 
 \paren{\nabla_{\mathbf{q}} \density{\specvec}(\mathbf{q},t) + \density{\specvec}(\mathbf{q},t)\nabla_{\mathbf{q}}\Psi^{\mathbf{s}}(\mathbf{q})}\right]
\end{equation*} 
for $L^{\mathbf{s}}$, where $I^n$ denotes the $n$ by $n$ identity matrix. For notational simplicity, the diffusion constants $D^{S_m} > 0$ have been factored out to multiply both the drift and diffusion terms (i.e., the potentials are dimensionless). The potential function $\Psi^{\mathbf{s}}$ may include a one-body potential $\psi^{S_m}$ for each species $S_m$, and interaction potentials $\psi^{S_mS_n}$ between each pair of particle species $S_m,S_n$: \begin{align}\label{eq:full_pot_fn}
    \Psi^{\mathbf{s}}(\mathbf{q}) =
 \sum_{m=1}^M\sum_{i=1}^{s_m} \left[ \psi^{S_m}(q_i^{s_m}) + \sum_{j=i+1}^{s_m} \psi^{S_mS_m}(q_i^{s_m},q_j^{s_m}) + \sum_{n < m} \sum_{j=1}^{s_n} \psi^{S_mS_n}(q_i^{s_m},q_j^{s_n})\right].
\end{align}

\subsection{Model problem}\label{section:introduction_reversible}

For simplicity of exposition, we present the new aspects of our CRDME using a representative model problem with three types of molecules ($M=3$), $S_1 = A$, $S_2 = B$, and $S_3=C$, a one-body potential $\psi^{S_m}$ for each species, an interaction potential $\psi^{S_mS_n}$ between each pair of species (including self-interactions), and the reversible reaction 
\begin{equation*}
A + B \rightleftharpoons C.
\end{equation*}
We emphasize, however, that our numerical approximation method applies to general reaction systems containing reversible and irreversible reactions of order $\leq 2$, with arbitrary combinations of one- and two-body potentials. Table~\ref{tab:discrete_reaction_rates} serves as a compact summary of the main results of our paper, the CRDME transition rates for spatial transport as well as a variety of zeroth- through second-order reactions, including the reversible $A + B \rightleftharpoons C$ reaction.

Let $\mathbf{s} = (a,b,c)$ denote the vector of copy numbers of the three species. To write the reaction PIDE operator $R^{\mathbf{s}}$, we define forward and backward reaction operators, respectively denoted by $\fwdop{\alpha}{\beta}{z}$ and $\bwdop{x}{y}{\xi}$, which act on the particle position vectors $\plocs = \slocsabc$. The forward operator acts on $\plocs$ by removing the $\alpha$th particle of species $A$ and the $\beta$th particle of species $B$, and adding a particle of species $C$ at position $z$ to the end of $\qc$. Similarly, $\bwdop{x}{y}{\xi}$ removes the $\xi$th particle of type $C$ and appends particles of type $A$ and $B$ at $x$ and $y$ to $\qa$ and $\qb$.

We will use $\fwdop{\alpha}{\beta}{z}$ and $\bwdop{x}{y}{\xi}$ in writing reactive interactions as functions of the pre-reaction and post-reaction position state vectors. For example, reaction terms for a \textit{forward} $A_\alpha + B_\beta \to C_z$ reaction from the starting state $\plocs$ will be a function of $(\fwdop{\alpha}{\beta}{z}\plocs,\plocs)$ (and will often also explicitly use the positions of the substrates, i.e., $(q^a_\alpha, q^b_\beta)$ here, and products, i.e., $z$ here, contained within these vectors). When there is no need to specify the locations and indices of the substrate and product particles, we will occasionally use the shorthand $\qplus$ and $\qminus$ in place of $\fwdop{\alpha}{\beta}{z}\plocs$ and $\bwdop{x}{y}{\xi}\plocs$. We also use the notation $\specvec^+ = (a-1,b-1,c+1)$ and $\specvec^- = (a+1,b+1,c-1)$ for the total particle counts in a pre- or post-reaction state. 

In this work, for reversible reactions we always factor the reactive interaction kernels into two pieces: a \textit{base reaction kernel} $\kappa^+(z|x,y)$ or $\kappa^-(x,y|z)$ depending only on the locations of the substrate and product particles, and a \textit{Boltzmann factor} $\nu^+(\fwdop{\alpha}{\beta}{z}\plocs, \mathbf{q})$ or $\nu^-(\bwdop{x}{y}{\xi}\plocs,\plocs)$ which takes the positions of all particles in the post-reaction and pre-reaction states as its two arguments. Typically, $\kappa^+(z|x,y)$ and $\kappa^-(x,y|z)$ represent the probability density per time for a binding $A_x + B_y \to C_z$ or unbinding $C_z\to A_x + B_y$ reaction to occur in the absence of any potentials. The multiplicative corrections $\nu^+(\fwdop{\alpha}{\beta}{z}\plocs, \mathbf{q})$ and $\nu^-(\bwdop{x}{y}{\xi}\plocs,\plocs)$ to the base reaction kernels then enforce detailed balance of reactive fluxes in the presence of potentials (see Section~\ref{section:continuum_detailed_balance}), and in particular are not required for irreversible reactions.

In our examples, we will ultimately follow a similar formulation to~\cite{frohner_reversible_2018,Heldman_2023} and choose the Boltzmann factor such that the reaction process can be interpreted as involving a probability density per time that a set of substrate particles will be proposed to react to produce a set of product particles, and a probability that the proposed reaction will be accepted. 

With the notations described above in mind, the reaction operator $\Rabc = R^{\mathbf{s}}$ for the model problem is given by 
\begin{equation}\label{eq:reversible-operator}
    \begin{split}
        (R^{\mathbf{s}}\vecdensity)(\vqa,\vqb,\vqc,t) &= - \density{\specvec}(\mathbf{q},t)\sum_{\alpha=1}^a\sum_{\beta=1}^b \int_\Omega \kappa^+(z|q^a_\alpha,q^b_\beta)\nu^+\left(\fwdop{\alpha}{\beta}{z}\plocs, \mathbf{q}\right)dz 
        \\&\emptytext{=} -\density{\specvec}(\plocs,t) \sum_{\xi=1}^c \int_{\Omega^2}\kappa^-(x,y|q^c_\xi)\nu^-\left(\bwdop{x}{y}{\xi}\plocs,\plocs\right)dx\,dy
        \\&\emptytext{=}+  \sum_{\xi=1}^c\int_{\Omega^2} \kappa^+(q^c_\xi | x,y)\nu^+(\mathbf{q}, \bwdop{x}{y}{\xi}\plocs )\density{\specvec^-}(\bwdop{x}{y}{\xi}\plocs,t)\,dx\,dy
        \\&\emptytext{=}+ \sum_{\alpha=1}^a\sum_{\beta=1}^b \int_\Omega \kappa^-(q^a_\alpha,q^b_\beta | z)\nu^-(\mathbf{q},\fwdop{\alpha}{\beta}{z}\plocs)\density{\specvec^+}\left(\fwdop{\alpha}{\beta}{z}\plocs,t\right)dz.
    \end{split}
\end{equation} We present the discretization of \eqref{eq:reversible-operator} in Section~\ref{section:derivations_reactions}.

\subsubsection{Continuum detailed balance conditions}\label{section:continuum_detailed_balance}
Since $\Omega$ is a closed domain and we are assuming reflecting Neumann boundary conditions, at equilibrium we expect the principle of detailed balance to hold for the reversible $A + B \rightleftharpoons C$ reaction~\cite{frohner_reversible_2018,zhang_detailed_2022}. This means that at equilibrium, the pointwise probability flux of forward reactions $A + B \to C$ should be exactly balanced by the pointwise probability flux of backward reactions $C \to A + B$. The detailed balance condition for transitions between the state $\mathbf{q}$, $\mathbf{s} = (a,b,c)$ and $\qplus = \fwdop{\alpha}{\beta}{z}\mathbf{q}$, $\mathbf{s}^+ = (a-1,b-1,c+1)$ is \begin{equation}\label{eq:DB-multiparticle-reversible}
\denseq(\mathbf{q})\kappa^+(z | {q}^a_\alpha,{q}^b_\beta)\nu^+(\mathbf{q}^+,\mathbf{q}) = \denseqplus(\mathbf{q}^+)\kappa^-( {q}^a_\alpha,{q}^b_\beta | z)\nu^-(\mathbf{q},\mathbf{q}^+),
\end{equation}
where we use $\denseq(\mathbf{q})$ to denote the equilibrium density. In words,
\eqref{eq:DB-multiparticle-reversible} means that at equilibrium, the
probability density per time of a forward reaction causing the transition from state
$\mathbf{q}$ to $\mathbf{q}^+$ is exactly balanced by the probability density per time of
the backward transition from $\mathbf{q}^+$ to $\mathbf{q}$. 

Since \eqref{eq:DB-multiparticle-reversible} implies $R^{\specvec} \bar{\vecdensity} = 0$ for all $\specvec$ and the time derivative in~\eqref{eq:FKE-full} vanishes at steady state, we find that $L^{\mathbf{s}}\denseq = 0$ as well. Let $\Psi^{\mathbf{s}}$ be the potential from \eqref{eq:full_pot_fn}. Assuming that $\Psi^{\bs}$ is sufficiently regular, the kernel of $L^{\mathbf{s}}$ is the one-dimensional space spanned by the Gibbs-Boltzmann density, $e^{-\Psi^{\bs}(\vq)}$. Taking into account our normalization convention, the equilibrium density is given by
\begin{equation*}
\denseq(\mathbf{q}) = \mathbf{s}!\frac{\bar{\pi}_\mathbf{s}}{Z_{\mathbf{s}}}e^{-\Psi^{\mathbf{s}}(\mathbf{q})},
\end{equation*}
where $\bar{\pi}_{\mathbf{s}}$ is the equilibrium value of \eqref{eq:particle_count_probabilities}, and 
\begin{equation*}
    Z_{\mathbf{s}} = \int_{\Omega^{|\mathbf{s}|}} e^{-\Psi^{\mathbf{s}}(\mathbf{q})}d\mathbf{q}
\end{equation*}
is the partition function. 

 Rearranging \eqref{eq:DB-multiparticle-reversible} and assuming, for an equilibrium constant $\Kd > 0$, the proportionality relation 
\begin{equation*}
\kappa^+(z|x,y) = \Kd \kappa^-(x,y|z),
\end{equation*}
as required for $\kappa^+$ and $\kappa^-$ to be in detailed balance in the purely diffusive case~\cite{zhang_detailed_2022}, we obtain  
\begin{equation} \label{eq:DB-full-equation}
    \frac{\nu^-(\mathbf{q},\mathbf{q}^+)}{\nu^+(\mathbf{q}^+,\mathbf{q})} = \Kd\frac{ab}{c + 1}\frac{\bar{\pi}_{\mathbf{s}}Z_{\mathbf{s}^+}}{\bar{\pi}_{\mathbf{s}^+}Z_{\mathbf{s}}}e^{-[\Psi^{\mathbf{s}}(\mathbf{q}) - \Psi^{\mathbf{s}^+}(\mathbf{q}^+)]}.
\end{equation}
Hence, for $\mathbf{s}$ fixed, as a function of $\mathbf{q}$ and $\mathbf{q}^+$ we have that 
\begin{equation}\label{eq:DB-ratio}
     \frac{\nu^-(\mathbf{q},\mathbf{q}^+)}{\nu^+(\mathbf{q}^+,\mathbf{q})} \propto e^{-[\Psi^{\mathbf{s}}(\mathbf{q}) - \Psi^{\mathbf{s}^+}(\mathbf{q}^+)]}.
\end{equation}

From \eqref{eq:DB-ratio} one sees that the more (less) potential state $\mathbf{q}^+$ has relative to state $\mathbf{q}$, the less (more) likely the forward reaction becomes relative to the backward reaction between the same two states, so that gradients in potential energy between post- and pre-reaction states modulate the relative likelihoods of the two reactions.

We discuss specific choices of $\nu^+$ and $\nu^-$ in Section~\ref{section:simulation_algorithm_pair_potentials_reactions}.

%% file: manuscript/background.tex
\section{CRDME background: one-body drift, diffusion, and pair reaction rates}\label{section:simulation_algorithm_background}
In this section we discuss approaches for deriving CRDME approximations to the PBSRD model \eqref{eq:FKE-full} on unstructured meshes in the absence of two-body potential interactions. We begin by describing the mesh structure and notation we will use throughout the paper, and then review methods for discretizing the drift-diffusion operator $L^{\mathbf{s}}$ in Section~\ref{section:simulation_algorithm_background_drift_diffusion} and the reaction operator $R^{\mathbf{s}}$ in Section~\ref{section:simulation_algorithm_background_reactions}.

Given a division $\mathcal T_h$ of the domain $\Omega\subseteq \RR^d$ into triangles ($d=2$) or tetrahedra ($d=3$), called the primal elements of the mesh, we associate to each vertex $x_i \in \Omega_h = \{x_j\}_{j=1}^N$ a polygonal/polyhedral control volume $V_i$ such that $\Omega = \cup_{i=1}^N V_i$. We call the resulting subdivision $\{V_i\}_{i=1}^N$ the dual mesh, and the $V_i$ the voxels or dual mesh elements. We denote by $E_{i j}$ the directed edge from vertex $x_i$ to its neighbor $x_j$. 

A simple definition of the dual mesh elements on a standard uniform triangulation or tetrahedralization of a box domain, which we will call the Cartesian dual, uses equally sized box-shaped control volumes. Typical choices on unstructured triangulations include the barycentric, or median, dual mesh, in which the boundaries of the dual mesh elements are defined by connecting the midpoint of each primal mesh edge (and the barycenter of each face when $d=3$) to the barycenters of its supporting elements~\cite{engblom_simulation_2009}, and the Voronoi dual mesh formed by connecting the circumcenters of adjacent primal mesh elements \cite{yates_importance_2013}. An example of a barycentric dual mesh (blue, dashed) resulting from a triangulation (black) of a disk domain is shown on the left in Figure~\ref{fig:dual_mesh} below, along with a close-up of the primal and dual mesh elements supporting a single primal mesh node. 

For both the barycentric dual and the Cartesian dual, each $V_i$ contains a single vertex of the original triangulation. For the barycentric dual, moreover, $|V_i \cap T| = \frac{|T|}{d + 1}1_{x_i\in T}$ for each $T \in \mathcal{T}_h$, and for both the barycentric and Cartesian duals (for the latter, at nodes away from the domain boundary, hence everywhere on a periodic domain) the voxel volumes satisfy
\begin{align}\label{eq:voxel-identities}
    |V_i| = \sum_{\{T \in \mathcal{T}_h \,|\, x_i \in T\}} \frac{|T|}{d + 1}.
\end{align}
We subsequently assume our mesh is such that~\eqref{eq:voxel-identities} holds. This implies that 
\begin{equation*}
\int_{\Omega} u_h(x)dx = \sum_{i=1}^N |V_i|u_h(x_i)
\end{equation*}
whenever $u_h$ is piecewise linear on $\mathcal{T}_h$, making it possible to write quadrature rules on $\Omega$ in terms of $d$-dimensional voxel measures. 

\begin{figure}[h]
    \centering
    \includegraphics[scale=.2]{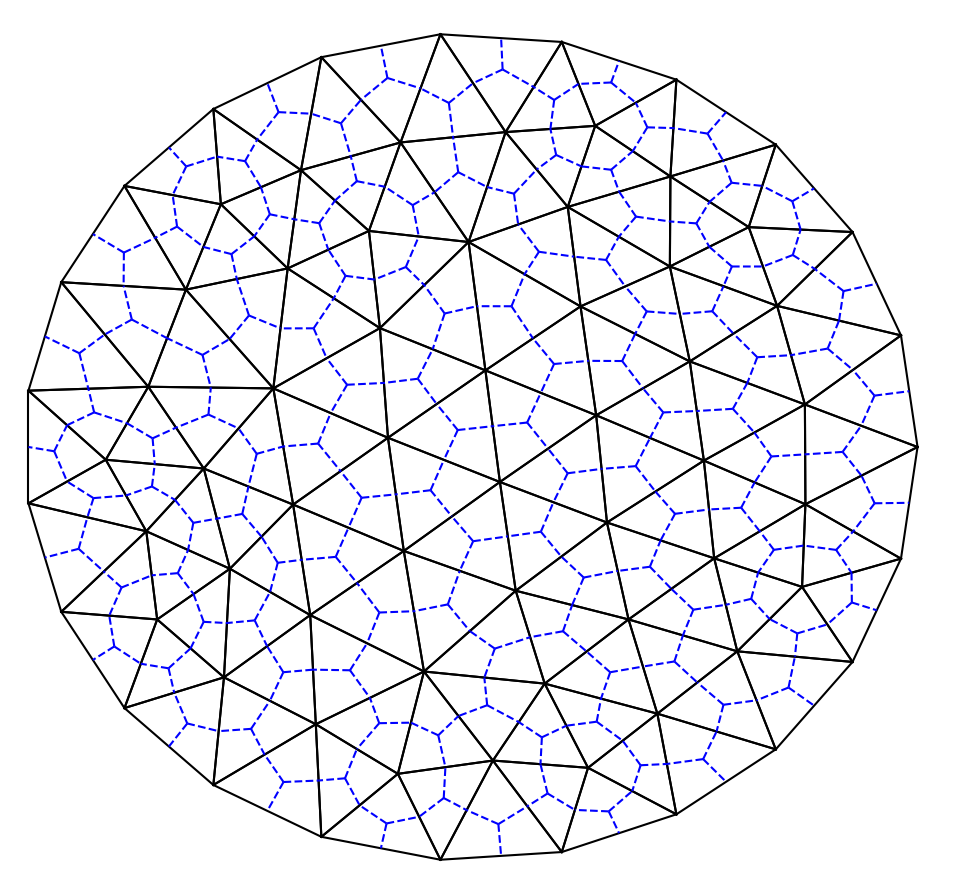}\hspace{10pt}\includegraphics[scale=.25]{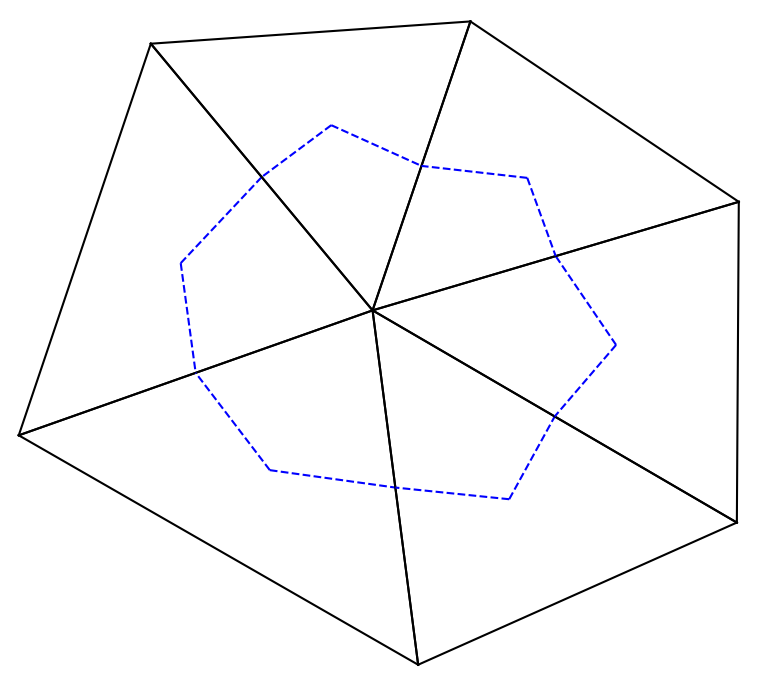}
    \caption{Left: A finite element primal mesh on a disk (black) and its barycentric dual (blue, dashed). Right: The supporting primal and dual mesh elements containing an interior mesh node.}
    \label{fig:dual_mesh}
\end{figure}

The jump process corresponding to our spatial discretization of the PBSRD model~\eqref{eq:FKE-full} will start with an initial distribution of particles of each type among the voxel sites. During a simulation, the particles move by discrete hops between neighboring dual mesh voxels, or equivalently along edges of the primal mesh, and reactions between particles in nearby voxels occur at certain rates. The hopping and reactive transition rates may depend on the distribution of particles in nearby voxels, in addition to the primal and dual mesh geometries.

The CRDME identifies hopping and reactive transition rates in such a way that statistics of the simulation (e.g., the expected number of particles of type $S_i$ at time $t$) converge to their continuum analogues (as described by the density solving \eqref{eq:FKE-full}) as the primal mesh is refined. Analogously to $\bq$, let $\multinode{i}{s} = (\bi^{s_1},\dots,\bi^{s_M}) \in \{1,\dots,N\}^{|\bs|}$ label the vector of mesh voxels in which a set of particles are located. Starting with the weak form of the PIDE \eqref{eq:FKE-full}, we discretize in space using a finite element method, yielding a coupled system of ODEs for $\cdfvec_h(t) = [\cdf{\specvec}_{\multinode{i}{s}}(t)]_{\specvec \in \NN_0^M}$. Here $\cdf{\specvec}_{\multinode{i}{s}}(t) : \{1,\dots,N\}^{|\bs|} \times [0,\infty) \to [0,\infty)$, with $\cdf{\specvec}_{\multinode{i}{s}}(t) / \bs! \in [0,1]$ representing the probability of observing $\specvec$ particles located at the lattice sites $\multinode{i}{s}$ at time $t$ (see Section~\ref{section:derivations} for details). Our final semi-discretization of~\eqref{eq:FKE-full} will then have the form
\begin{align}\label{eq:discrete-FKE-full}
 \derivt{\cdf{\specvec}_{\multinode{i}{s}}}(t) = [L^{\mathbf{s}}_h \cdf{\specvec}]_{\multinode{i}{s}}(t) + [R^{\mathbf{s}}_h \cdfvec_h]_{\multinode{i}{s}}(t),
\end{align} 
representing a coupled system of ODEs over all possible number states, $\bs$, and all lattice particle configurations for each number state, $\multinode{i}{s}$. This discretization can be interpreted as the forward equation for a spatial jump process, i.e., a continuous-time Markov chain (CTMC) on the mesh, under certain mesh conditions which we discuss in more detail later in this section. When this interpretation is valid, the associated jump process can be exactly sampled using the stochastic simulation algorithm (SSA). The hopping transition rates between voxels are contained in the off-diagonal elements of $L^{\mathbf{s}}_h$, and the reactive transition rates in $R^{\mathbf{s}}_h$.
 
\subsection{One-body drift and diffusion}\label{section:simulation_algorithm_background_drift_diffusion}

When the system does not include two-body interaction potentials, it is common to approximate the transport and reaction operators separately. Since in this case the transport involves no particle interactions, particles move independently (in the absence of reactions), and it is sufficient to derive spatial hopping transition rates for the transport of one particle~\cite{IsaacsonSisc2006}. 
We therefore start with the drift-diffusion equation for the probability density, $f^S(x,t)$, that a particle of species $S$ is at $x \in \Omega$ at $t$, 
\begin{equation}\label{eq:FKE-drift-diffusion-single}
\partialt{\densityletter^{S}} = D^S \nabla \cdot (\nabla \densityletter^{S} + \densityletter^S \nabla \psi^S).
\end{equation} 
One then discretizes this equation to a master equation for a spatial jump process, for which the individual hopping rates for particles of that species correspond to the discretization coefficients. See~\cite{engblom_simulation_2009,isaacson_convergent_2013,isaacson_unstructured_2018} for the pure diffusion case where $\psi^S$ is constant, and the more recent paper \cite{isaacson_zhang_crddme_2024} when including drift due to one-body potentials. 

In the multiparticle system, this procedure can be interpreted as a discretization of the operator $L^{\mathbf{s}}$ from \eqref{eq:FKE-full} by applying a tensor product FEM to \eqref{eq:FKE-full}, and then repeatedly applying mass-lumping quadrature \eqref{eq:mass-lump-quadrature} to ensure that the discrete operator $L^{\mathbf{s}}_h$ decouples into a sum of single-particle hopping operators (a Kronecker sum), one for each particle. We expand on these ideas, which have not appeared in the literature, in Section~\ref{section:simulation_algorithm_pair_potentials}, where the inclusion of two-body interaction potentials implies that the motion of an individual particle cannot be fully decoupled from other particles. That is, the hopping rates for a given particle will now depend on the locations of other particles in the system.

Dropping the species labels for now in \eqref{eq:FKE-drift-diffusion-single}, in the purely diffusive case that $\psi$ is constant, for a single particle one simply applies a $P^1$-FEM discretization to \eqref{eq:FKE-drift-diffusion-single} to obtain a discrete Laplacian $\Delta_h = DW_h (M_h^L)^{-1}$, where $W_h$ is the Laplacian stiffness matrix defined below and $M_h^L$ is the lumped mass matrix $(M_h^L)_{ij} = |V_i|\delta_{ij}$, where 
\begin{equation*}
\delta_{ij} = 
    \begin{cases} 1 & i = j\\
        0 & \text{otherwise}
    \end{cases}.
\end{equation*}
The hopping rate from voxel $V_i$ to voxel $V_j$ is given by $(\Delta_h)_{ji} = D\frac{\omega_{ji}}{|V_i|}$, where
\begin{equation*}
\omega_{ij} = (W_h)_{ij} = -\int_\Omega \nabla\phi_i \cdot \nabla\phi_j \, dx,
\end{equation*}
with $\{\phi_i\}_{i=1}^N$ the $P^1$ nodal basis on $\mathcal{T}_h$, as first derived in constructing an unstructured-grid RDME in~\cite{engblom_simulation_2009}. In the remainder, we use $\omega_{ij}$ to denote the entries of this matrix, i.e., the negative of the standard $P^1$ stiffness matrix for the Laplacian, so that $\omega_{ij} \geq 0$ for $i \neq j$ on suitable meshes and $\sum_{i=1}^N \omega_{ij} = 0$.

For the hopping rates to be valid, they must be nonnegative, i.e., $\omega_{ij}\geq 0$ for $i \neq j$. Generally, for $\Delta_h$ to represent the transition rate matrix of a CTMC it must also have zero column sums. Hence, we need $(\Delta_h)_{ij} \geq 0$ for $i\neq j$ and 
\begin{equation*}
\sum_{i=1}^N (\Delta_h)_{ij} = \sum_{i=1}^N D\frac{\omega_{ij}}{|V_j|} = \frac{D}{|V_j|} \sum_{i=1}^N \omega_{ij} = 0,\quad j=1,\dots,N.
\end{equation*}
The latter property is satisfied for typical FEM discretizations, for which the sum of the basis functions is a constant. One can also show that the Laplacian stiffness matrix $W_h$ has nonnegative off-diagonals whenever the mesh satisfies a particular condition on the dihedral angles of $\mathcal{T}_h$, which, e.g., Delaunay triangulations satisfy for $d=2$ \cite{xu_monotone_1999,engblom_simulation_2009}.

On the other hand, if $\psi$ is nonconstant then the analogous Galerkin discretization of the drift-diffusion operator does not necessarily result in nonnegative hopping rates, even if the triangulation is Delaunay. For the particular case of a one-body potential, the problem can be circumvented by using a specific quadrature formula, the EAFE method described in \cite{xu_monotone_1999}, to generate a stiffness matrix with entries 
\begin{equation}\label{eq:stiffness_background_pot} (A_h)_{ij} = D\omega_{ij} \left[\frac{1}{|E_{ij}|}\int_{E_{ij}} e^{\psi(x) - \psi(x_j)}dx\right]^{-1}
\end{equation} 
if $i\neq j$, and $(A_h)_{ii} = -\sum_{j \neq i} (A_h)_{j i}$. Evidently the matrix $A_h$ has nonnegative off-diagonals and zero column sums whenever the Laplacian $(\omega_{ij})_{i,j=1}^N$ does \cite{xu_monotone_1999}. This approach was used to derive a CRDME discretization of~\eqref{eq:FKE-full} in the case of reaction-drift-diffusion systems with only one-body potentials~\cite{isaacson_zhang_crddme_2024}.

The EAFE method as derived in \cite{xu_monotone_1999} is equivalent to the finite element method described in \cite{markowich_inverse-average-type_1988} for the purpose of semiconductor device modeling, and can be interpreted as a finite element generalization of the finite difference inverse-average scheme for 1D semiconductor modeling first appearing in \cite{scharfetter_large-signal_1969}. Indeed, to obtain an exponentially weighted inverse-average scheme in the spirit of the latter, one can replace $\psi$ by its piecewise linear nodal interpolant in \eqref{eq:stiffness_background_pot} to get the formula 
\begin{equation}\label{eq:stiffness_background_pot_quad}
    (A_h)_{ij} = D\omega_{ij} \frac{\bd_{ji}\psi}{e^{\bd_{ji}\psi} - 1} = D\omega_{ij}B(\bd_{ji}\psi), \quad i \neq j,
\end{equation} 
where we use the difference operator notation $\bd_{ji} \psi := \psi(x_i) - \psi(x_j)$ 
and 
\begin{equation}\label{eq:bernoulli}
    B(z) =\begin{cases} \frac{z}{e^z - 1} & z \neq 0\\
                            1 & z = 0 \end{cases} 
\end{equation}
is the Bernoulli function. For mesh functions $u_h(x_i)$ and $v_h(x_i)$, the corresponding discretization of the bilinear form in the weak representation of~\eqref{eq:FKE-drift-diffusion-single} is then
\begin{equation}\label{eq:EAFE_bilinear_form}
    a_h(u_h,v_h) = -\sum_{E_{ij} \in \mathcal{T}_h} D\omega_{ij}B(\bd_{ij}\psi)\bd_{ij}\left(e^{\psi - \psi(x_i)}u_h\right)\bd_{ij}v_h,
\end{equation} 
where the sum is over all distinct edges $E_{ij}$ in the triangulation $\mathcal{T}_h$, and each edge is given one (arbitrary) orientation and appears only once within the sum. The modified method \eqref{eq:stiffness_background_pot_quad} also appears in \cite{markowich_inverse-average-type_1988}, who mention that in applications the potential $\psi$ often comes from the $P^1$-finite element solution of a Poisson equation. We always use analogues of \eqref{eq:stiffness_background_pot_quad} (in place of \eqref{eq:stiffness_background_pot}) and the resulting bilinear form \eqref{eq:EAFE_bilinear_form} in our discretizations to avoid computing the high-dimensional integrals that would result from using \eqref{eq:stiffness_background_pot} with multiparticle equations involving interaction potentials (see Section~\ref{section:derivations_hopping}). %

Letting $\bar{u}(x) = e^{-\psi(x)}$, the Gibbs-Boltzmann distribution associated with $\psi$, we observe that 
$\bar{u}$ is the equilibrium solution to 
\eqref{eq:FKE-drift-diffusion-single}. %
Direct substitution shows that the EAFE discrete solution satisfies $\bar{u}_h(x_i) = \bar{u}(x_i)$, i.e., the EAFE scheme applied to~\eqref{eq:FKE-drift-diffusion-single} gives the exact steady-state nodal values. Moreover, because the steady-state solution only depends on the nodal values of $\psi$, if we use a quadrature formula for the edge integrals as in \eqref{eq:stiffness_background_pot_quad}, the steady-state solution is preserved. 

Hence, we summarize this discussion by concluding that whether we use the formula \eqref{eq:stiffness_background_pot} or \eqref{eq:stiffness_background_pot_quad} for the stiffness matrix entries associated with a discretization of \eqref{eq:FKE-drift-diffusion-single}, we obtain the following properties:

\begin{enumerate}
\item The stiffness matrix $A_h$ generates a transition rate matrix for a spatial jump process if the corresponding $P^1$-finite element Laplacian does.
\item The steady-state solution to the spatially discrete version of \eqref{eq:FKE-drift-diffusion-single} is $\bar{u}_i \propto e^{-\psi(x_i)}$, where $\bar{u}_{i} = \bar{u}_h(x_i)$ is the $i$th entry of the vector of coefficients for the EAFE finite element solution $\bar{u}_h$.
\end{enumerate}

For the interaction potential case, we will construct our approximation of the spatial transport operator such that these two properties also hold.

\subsection{Reaction rates}\label{section:simulation_algorithm_background_reactions}

In \cite{isaacson_unstructured_2018}, the authors use a finite volume method to discretize the reaction operator for the reversible $A + B \xrightleftharpoons[\kappa^-(x,y|z)]{\,\kappa^+(z|x,y)\,} C$ reaction in the purely diffusive case. Motivated by that and other recent work~\cite{isaacson_convergent_2013,zhang_detailed_2022,isaacson_zhang_crddme_2024}, we assume a specific factorization of the reaction kernels as 
\begin{align*}
    \kappa^+(z|x,y) &= \kappa(x,y) m(z|x,y), & \Kd \kappa^-(x,y|z) &= \kappa^+(z|x,y),
\end{align*}
where $\Kd > 0$ represents an equilibrium constant for the reaction. In the purely diffusive case, $\kappa^- \propto \kappa^+$ is required for detailed balance of pointwise reactive fluxes to hold at equilibrium~\cite{zhang_detailed_2022}, and $\Kd$ is typically chosen as the inverse of the well-mixed dissociation constant to ensure the equilibrium solution to the PBSRD model is consistent with that of the well-mixed chemical master equation~\cite{isaacson_unstructured_2018}. 

$\kappa(x,y)$ can be interpreted as the probability per time an $A$ at $x$ reacts with a $B$ at $y$, and $m(z|x,y)$ as the probability density for placing the resulting $C$ at $z$ given that the substrates were at $x$ and $y$ (in the purely diffusive case). $\kappa^-$ can be decomposed analogously. A common choice for $\kappa$ is the Doi reaction kernel \cite{doi_stochastic_1976}, given by
\begin{equation}\label{eq:doi_kernel}
    \kappa(x,y) = \lambda 1_{|x-y| < \varepsilon},
\end{equation} where $\lambda > 0$ is the microscopic probability per time of the reaction occurring when the two substrates are within a reaction radius $\varepsilon > 0$. There are many options for the placement kernel, but one common choice is 
\begin{equation}\label{eq:nodal_placement}
    m(z|x,y) = \left[\gamma \delta(x - z) + (1-\gamma)\delta(y-z)\right], 
\end{equation}
where with probability $\gamma$ the product $C$ molecule is placed at the location of the $A$ molecule, and with probability $1-\gamma$ it is placed at the location of the $B$ molecule. 
We note that the product placement mechanism \eqref{eq:nodal_placement} is similar to the one used in the Smoldyn PBSRD simulation software \cite{andrews_stochastic_2004,andrews_detailed_2010,andrews_smoldyn_2017}, but that other choices are possible (for example, placing at the diffusion-weighted center of mass~\cite{isaacson_unstructured_2018}). 

With the finite-volume-inspired discretization from \cite{isaacson_unstructured_2018} in the context of developing a spatial SSA, the discrete forward and backward reaction rates between an $A$ at $V_i$, a $B$ at $V_j$, and a $C$ at $V_k$ are
\begin{align}
\kappa^+_{ijk} &= \frac{1}{|V_{ij}|}
\int_{V_{ijk}}\kappa^+(z|x,y)dxdydz = \left[\gamma\delta_{ik} + (1-\gamma)\delta_{jk}\right]\frac{1}{|V_{ij}|}\int_{V_{ij}} \kappa(x,y)dxdy,\label{eq:discrete_fwd_standard}\\
\kappa^-_{kij} &= \frac{1}{|V_{k}|}\int_{V_{kij}} \kappa^-(x,y|z)dzdxdy = \frac{1}{\Kd^h}[\gamma \delta_{ik} + (1-\gamma)\delta_{jk}]\frac{1}{|V_k|}\int_{V_{ij}}\kappa(x,y)dxdy,\label{eq:discrete_bwd_standard}
\end{align} 
where $V_{ij} = V_i\times V_j$ and $V_{ijk} = V_i\times V_j\times V_k$ (with permuted subscripts denoting the correspondingly ordered product), so that $|V_{ij}| = |V_i||V_j|$ and $|V_{ijk}| = |V_i||V_j||V_k|$. Here we have introduced a mesh-dependent equilibrium constant $\Kd^h$ to be consistent with our later equations in the presence of interaction potentials. In the present, purely diffusive setting, it is simply the continuum equilibrium constant, i.e., $\Kd^h := \Kd$, but it will become distinct when we calibrate the discrete rates in the presence of interaction potentials. 
In Section~\ref{section:two_particle_equilibrium_constant} we show how to choose $\Kd^h$ in the presence of interaction potentials, in order to exactly reproduce a desired equilibrium ratio for the probability of being in the bound vs. unbound states in a two-particle system. 

For the simulations of Section~\ref{section:numerics_pair_potentials}, we use an approximation to these rates, obtained by applying a combination of the mass-lumping \eqref{eq:mass-lump-quadrature} and finite volume quadrature rules to these integrals, giving:
 \begin{align}
 \kappa^+_{ijk} &= \left[\gamma\delta_{ik}\frac{1}{|V_j|}\int_{V_j}\kappa(x_i,y)dy + (1-\gamma)\delta_{jk}\frac{1}{|V_i|}\int_{V_i}\kappa(x,x_j)dx\right],\label{eq:discrete_fwd_lump}\\
\kappa^-_{kij} &= \frac{1}{\Kd^h}\left[\gamma \delta_{ik} \int_{V_j} \kappa(x_k,y)dy + (1-\gamma) \delta_{jk}\int_{V_i}\kappa(x,x_k)dx\right].\label{eq:discrete_bwd_lump}
\end{align} 

Finally, we note that one may consider an alternative approximation to the rates by using a pure mass-lumping quadrature, to derive 
\begin{align}
    \kappa^+_{ijk} &= \kappa(x_i,x_j) [\gamma\delta_{ik} + (1-\gamma)\delta_{jk}],\label{eq:discrete_fwd_smooth}\\
    \kappa^-_{kij} &=  \frac{\kappa(x_i,x_j)}{\Kd^h}[\gamma\delta_{ik}|V_j| + (1-\gamma)\delta_{jk}|V_i|].\label{eq:discrete_bwd_smooth}
\end{align} 
For the discontinuous Doi reaction kernel, $O(h^2)$ convergence of simulation statistics has been observed with \eqref{eq:discrete_fwd_standard} and \eqref{eq:discrete_bwd_standard} in the drift-free case~\cite{isaacson_unstructured_2018}, and with \eqref{eq:discrete_fwd_lump} and \eqref{eq:discrete_bwd_lump} in Section~\ref{section:numerics_pair_potentials} and in~\cite{Heldman_2023}. In contrast, one expects that the pure mass-lumping rates \eqref{eq:discrete_fwd_smooth} and \eqref{eq:discrete_bwd_smooth} are only first-order accurate for discontinuous kernels; see~\cite{Heldman_2023} for additional discussion. For smooth reaction kernels, all three choices are expected to be second-order accurate.

Whether one uses the rates \eqref{eq:discrete_fwd_standard} and \eqref{eq:discrete_bwd_standard}, \eqref{eq:discrete_fwd_lump} and \eqref{eq:discrete_bwd_lump}, or \eqref{eq:discrete_fwd_smooth} and \eqref{eq:discrete_bwd_smooth}, the discrete rates satisfy the (pure diffusion case) detailed balance condition~\cite{isaacson_unstructured_2018}
\begin{align}\label{eq:discrete-DB-condition-no-potentials}
    |V_{ij}|\kappa^+_{ijk} = \Kd^h|V_k|\kappa^-_{kij}.
\end{align}
Building from the preceding discretized rates, in Section~\ref{section:simulation_algorithm_pair_potentials_reactions} we show how to generalize the PBSRD model reaction operator discretization such that detailed balance is preserved in the presence of interaction potentials.

\section{Hopping and reaction rates in the presence of two-body interaction potentials}\label{section:simulation_algorithm_pair_potentials}
\indent In this section we present our main results: the new CRDME, and associated hopping and reactive transition rates. They are obtained by discretizing the full PBSRD model with two-body potential interactions~\eqref{eq:FKE-full}. Spatial hopping transition rates are presented in Section~\ref{section:simulation_algorithm_pair_potentials_hopping_rates}, based on the derivation in Section~\ref{section:derivations_hopping}. We find that the derived spatial transition rates between two states now depend on the potential difference between them, reflecting that in the presence of two-body potentials the motion of individual particles may be impacted by many other particles in the system. The more detailed Section~\ref{section:derivations_hopping} also establishes that the derived rates are consistent with detailed balance of drift-diffusion fluxes holding at equilibrium, and that the CRDME supports an equilibrium Gibbs-Boltzmann distribution for the particle positions. In Section~\ref{section:simulation_algorithm_pair_potentials_reactions}, we present our discrete reactive transition rates, based on the derivation in Section~\ref{section:derivations_reactions}, and show how they satisfy a discrete reactive detailed balance condition consistent with that of the continuum PBSRD model. Section~\ref{section:two_particle_equilibrium_constant} then investigates how to choose a mesh-dependent equilibrium constant $\Kd^h$ to ensure consistency with the probability of being in the bound vs. unbound states for a two-particle system at equilibrium, motivated by the parametrization approach used for deriving a Brownian dynamics method in~\cite{frohner_reversible_2018}. Finally, we conclude in Section~\ref{section:summary} by summarizing the spatial hopping and reactive transition rates needed to define and simulate the spatial jump process associated with the new CRDME for general reaction systems (see Table~\ref{tab:discrete_reaction_rates}). Full derivations of the general multiparticle spatial hopping transition rates, the reactive transition rates, and a discussion of how the CRDME is consistent with the principle of detailed balance and supports Gibbs-Boltzmann spatial distributions at equilibrium are given in Section~\ref{section:derivations}. 

\subsection{Hopping rates}\label{section:simulation_algorithm_pair_potentials_hopping_rates}

Let $i^{s_m}_\alpha \in \{1,\dots,N\}$ denote the index of the voxel $V_{i^{s_m}_\alpha}$ containing particle $\alpha \in \{1,\dots,s_m\}$ of type $S_m$, and let $\mathbf{i}^{s_m} = (i_1^{s_m},\dots,i_{s_m}^{s_m})$ be the position state vector for particles of type $S_m$. The (mesh) position state vector for all particles when the number state of the system is $\mathbf{s}$ is then $\mathbf{i}^{\mathbf{s}} = (\mathbf{i}^{s_1},\dots,\mathbf{i}^{s_M})$, as in Section~\ref{section:simulation_algorithm_background}. Let $\Psi_{\mathbf{i}^{\mathbf{s}}} := \Psi^{\mathbf{s}}(\mathbf{x}_{\mathbf{i}^{\mathbf{s}}})$ be the value of $\Psi^{\mathbf{s}}$ from \eqref{eq:full_pot_fn} evaluated at the corresponding node $\mathbf{x}_{\mathbf{i}^{\mathbf{s}}} \in \Omega^{|\mathbf{s}|}$, where
\begin{equation*}
\mathbf{x}_{\mathbf{i}^{\mathbf{s}}} = (\mathbf{x}_{\mathbf{i}^{s_1}},\dots,\mathbf{x}_{\mathbf{i}^{s_M}})\qquad\text{and}\qquad \mathbf{x}_{\mathbf{i}^{s_m}} = (x_{i^{s_m}_1},\dots,x_{i^{s_m}_{s_m}}),\quad m=1,\dots,M,
\end{equation*}
for $\{x_i\}_{i=1}^N$ the nodes of the mesh. We use the notation \begin{equation}
    \hopop{s_m}{\alpha}{j}\mathbf{i}^{\mathbf{s}} = (\mathbf{i}^{s_1},\dots,\mathbf{i}^{s_{m-1}},i^{s_m}_1,\dots,i^{s_m}_{\alpha-1},j,i^{s_m}_{\alpha+1},\dots,i^{s_m}_{s_m},\mathbf{i}^{s_{m+1}},\dots,\mathbf{i}^{s_M})
\end{equation}
to represent the state of the system after a spatial hop of the $\alpha$th particle of species $S_m$ from voxel $V_{i^{s_m}_\alpha}$ to a neighboring voxel $V_j$.  Then, the transition rate for a specific particle with $i_{\alpha}^{s_m} = i$ to hop to voxel $V_j$ is given by 
\begin{align}\label{eq:discrete_hopping_rate_formula}
H^{S_m}_{ji}\left(\mathbf{x}_{\mathbf{i}^{\mathbf{s}}}\right) = D^{S_m}\frac{\omega_{j i}}{|V_i|}B\left(\Psi_{\hopop{s_m}{\alpha}{j}\mathbf{i}^{\mathbf{s}}} - \Psi_{\mathbf{i}^{\mathbf{s}}}\right),
\end{align}
i.e., the hopping rate depends, via the Bernoulli function \eqref{eq:bernoulli}, on the potential difference between the pre- and post-hop states. We present a detailed derivation of \eqref{eq:discrete_hopping_rate_formula} in Section~\ref{section:derivations_hopping} for a simplified two-species model. Since particles of the same species at the same position are indistinguishable, the total transition rate (i.e., propensity) for a particle of type $S_m$ in voxel $V_i$ to hop to voxel $V_j$ is simply 
\begin{equation*}
    H^{S_m}_{ji}(\mathbf{x}_{\mathbf{i}^{\mathbf{s}}}) s_{m i},
\end{equation*}
where $s_{m i}$ is the number of particles of type $S_m$ in voxel $V_i$. In the case that the potential is the same in both states, we recover the previously derived total diffusive transition rate, $D^{S_m}\frac{\omega_{j i}}{|V_i|} s_{m i}$~\cite{engblom_simulation_2009,isaacson_unstructured_2018}. 

In practice, during a simulation we keep track of the total potential difference $\Psi_{\hopop{s_m}{\alpha}{j}\mathbf{i}^{\mathbf{s}}} - \Psi_{\mathbf{i}^{\mathbf{s}}}$ from \eqref{eq:discrete_hopping_rate_formula} for each species and each directed mesh edge. This potential difference can be written as a running sum over the contributions from each particle in the system, so that when a particle is added to or removed from a voxel, we need only add or subtract the appropriate potential differences corresponding to the change in the system state. For each voxel, we also precompute and store a stencil of nearby voxels for which the stored potential differences need to be updated after an event; for the short-range potentials we consider in this work, the number of voxels that need to be updated after an event scales like $(\varepsilon_{\text{cutoff}}/h)^d$, where $\varepsilon_{\text{cutoff}}$ is the cutoff distance and $h$ is the mesh size. Since the accuracy of the discretization is governed by $\varepsilon_{\text{cutoff}}/h$ (see Section~\ref{section:numerics_pair_potentials}), this ratio can in practice be a moderate constant depending on the choice of potential function and desired accuracy, independent of the number of particles in the system.

\subsection{Reactions and detailed balance}\label{section:simulation_algorithm_pair_potentials_reactions}

In Section~\ref{section:continuum_detailed_balance}, we discussed how, with the inclusion of one-body or two-body interaction potentials, the base reaction kernels $\kappa^+(z|x,y)$ and $\kappa^-(x,y|z)$ are modulated by Boltzmann factors $\nu^+$ and $\nu^-$ depending on the potential energy of the system in the pre- and post-reaction states. In our discretized system, we likewise modulate the \textit{discrete} transition rates $\kappa^+_{ijk}$ and $\kappa^-_{kij}$ from Section~\ref{section:simulation_algorithm_background_reactions} with \textit{nodal} Boltzmann factors. 

In Section~\ref{section:derivations}, we derive discretized reactive transition rates when there are many particles of types $A$, $B$, and $C$. We again let $\multinode{i}{s} = (\mathbf{i}^a,\mathbf{i}^b, \mathbf{i}^c)$ be the particle positions, and let $\mathbf{x}_{\multinode{i}{s}}$ be the associated vector of mesh nodes. We also abuse notation and apply the operators $\fwdop{\alpha}{\beta}{z}$ and $\bwdop{x}{y}{\xi}$ in the obvious way to mesh nodes $\node{\multinode{i}{s}}$. Then, the derived reactive transition rates in the presence of potentials are given by 
\begin{equation}\label{eq:reaction-rates-general-form}
\begin{split}
    \hat{\kappa}^+_{i^a_\alpha i^b_\beta k}\left(\multinode{i}{s}\right) &= \nu^+\left(\fwdop{\alpha}{\beta}{k}\node{\multinode{i}{s}},\node{\multinode{i}{s}}\right)\kappa^+_{i^a_\alpha i^b_\beta k},\\
    \hat{\kappa}^-_{i^c_\xi ij}(\multinode{i}{s}) &= \nu^-\left(\bwdop{i}{j}{\xi}\node{\multinode{i}{s}},\node{\multinode{i}{s}}\right)\kappa^-_{i^c_\xi ij},
\end{split}
\end{equation}
with $\nu^+$, $\nu^-$ as in \eqref{eq:DB-ratio}. Here $\hat{\kappa}^+_{i^a_\alpha i^b_\beta k}\left(\multinode{i}{s}\right)$ can be interpreted as the probability per time one specific $A$ at $x_{i^a_\alpha}$ reacts with one specific $B$ at $x_{i^b_\beta}$ to form a $C$ at $x_k$, given that the system is in state $\multinode{i}{s}$. If there are $a_i$ particles of type $A$ in $V_i$ and $b_j$ particles of type $B$ in $V_j$, the total transition rate (i.e., propensity) for a reaction to occur between any $A$ in $V_i$ and any $B$ in $V_j$ to form a $C$ in $V_k$ is then $\hat{\kappa}^+_{i j k}\left(\multinode{i}{s}\right) a_i b_j$. $\hat{\kappa}^-_{i^c_\xi ij}(\multinode{i}{s})$ can be interpreted similarly. 

Similar to the reactive interaction kernels of the underlying PBSRD model, the discrete reactive transition rates of the new CRDME defined by~\eqref{eq:reaction-rates-general-form} are given by the product of the (base) reactive transition rates in the absence of potentials with the Boltzmann factors $\nu^+$ and $\nu^-$ evaluated at the nodal particle positions corresponding to the pre- and post-reaction states. 

\begin{remark}
The factorization \eqref{eq:reaction-rates-general-form} makes possible the use of a rejection method to simulate bimolecular reactions with Doi reaction kernels, instead of precomputing the integrals in \eqref{eq:discrete_fwd_standard}--\eqref{eq:discrete_bwd_lump}. See Chapter 3 of the thesis \cite{Heldman_2023} %
for details. 
\end{remark}

As concrete examples of $\nu^+$ and $\nu^-$, we define for each species pair $S_m,S_n$ the combined potential 
\begin{equation*}
\tilde{\psi}^{S_mS_n}(x,y) := \psi^{S_m}(x) + \psi^{S_n}(y) + \psi^{S_mS_n}(x,y),
\end{equation*}
set $\mathbf{q}^+ = \fwdop{\alpha}{\beta}{k}\qvec$ and $\mathbf{q}^- = \bwdop{i}{j}{\xi}\qvec$, and follow \cite{frohner_reversible_2018} by choosing 
\begin{equation}\label{eq:Noe-DB-boltzmann}
\begin{split}
    \nu^+(\mathbf{q}^+,\mathbf{q}) &=  e^{-\psi^C(z)}\min\left\{1, e^{-[\Psi^{\mathbf{s}^+}(\mathbf q^+) - \Psi^{\mathbf{s}}(\mathbf q)] + \psi^C(z) - \tilde{\psi}^{AB}(x_{i^a_\alpha},x_{i^b_\beta})}\right\},\\
    \nu^-(\mathbf{q}^-,\mathbf{q}) &= e^{-\tilde{\psi}^{AB}(x,y)}\min\left\{1, e^{-[\Psi^{\mathbf{s}^-}(\mathbf q^-) - \Psi^{\mathbf{s}}(\mathbf q)] + \tilde{\psi}^{AB}(x,y) - \psi^C(x_{i^c_\xi})}\right\}.
\end{split}
\end{equation}
Here we have used $x,y,z\in\RR^d$ to denote the positions of the three particles involved in the reaction, and slightly generalized \cite{frohner_reversible_2018} to include one-body potentials in the same manner as two-body potentials. In the remainder we will assume these specific choices for $\nu^+$ and $\nu^-$, the only choices that have been suggested in the literature to our knowledge, but note that other choices are possible as long as they are consistent with~\eqref{eq:DB-full-equation}.

Here and in Table~\ref{tab:discrete_reaction_rates}, we abbreviate $\Psi^{\bs}(\mathbf{q})$ as $\Psi(\mathbf{q})$, the number state being determined by $\mathbf{q}$. Once we insert \eqref{eq:Noe-DB-boltzmann} into \eqref{eq:reaction-rates-general-form}, the fact that \begin{align*}
\begin{array}{l}
   \alpha^+(\qplus,\qvec) := \min\left\{1, e^{-[\Psi(\mathbf q^+) - \Psi(\mathbf q)] + \psi^C(z) - \tilde{\psi}^{AB}(x_{i^a_\alpha},x_{i^b_\beta})}\right\} \\
    \alpha^-(\qminus,\qvec) := \min\left\{1, e^{-[\Psi(\mathbf q^-) - \Psi(\mathbf q)] + \tilde{\psi}^{AB}(x,y) - \psi^C(x_{i^c_\xi})}\right\}
\end{array}
\leq 1    
\end{align*} 
allows a new decomposition 
\begin{align*}
     \hat{\kappa}^+_{{i^a_\alpha}{i^b_\beta}k}({\multinode{i}{s}}) &= \tilde{\kappa}^+_{{i^a_\alpha}{i^b_\beta}k}\alpha^+\left(\fwdop{\alpha}{\beta}{k}\node{\multinode{i}{s}},\node{\multinode{i}{s}}\right),\\
    \hat{\kappa}^-_{i^c_\xi ij}({\multinode{i}{s}}) &= \tilde{\kappa}^-_{i^c_\xi ij}
     \alpha^-\left(\bwdop{i}{j}{\xi}\node{\multinode{i}{s}},\node{\multinode{i}{s}}\right),
\end{align*}
where 
\begin{equation}\label{eq:DB-rates-discrete-two-particle}
    \begin{split}
        \tilde{\kappa}^+_{ijk} &= e^{-\psi^C(x_k)}\kappa^+_{ijk},\\
        \tilde{\kappa}^-_{kij} &= e^{-\tilde{\psi}^{AB}(x_i,x_j)}\kappa^-_{kij}.
    \end{split}
    \end{equation} 
Here we can interpret $\tilde{\kappa}^+_{ijk}$ as representing a reaction proposal transition rate, i.e., the probability per time a reaction is proposed to occur between one specific $A$ in voxel $V_i$ and one specific $B$ in voxel $V_j$ to form a $C$ in voxel $V_k$, with $\tilde{\kappa}^-_{kij}$ defined analogously. The potential terms in~\eqref{eq:DB-rates-discrete-two-particle} modulate the purely diffusive reactive transition rates ($\kappa^+_{ijk}$ and $\kappa^-_{kij}$) to account for potentials induced by the product particles. The acceptance factors $\alpha^+$ and $\alpha^-$ account for the influence of particles beyond just the substrates and products. Within a CRDME simulation, we then use $\tilde{\kappa}^+_{ijk}$ and $\tilde{\kappa}^-_{kij}$ as reaction \textit{proposal} rates for reversible reactions and $\alpha^+$ and $\alpha^-$ as \textit{acceptance} probabilities for proposed reactions. We explain this concept further in our results summary in Section~\ref{section:summary}.

\subsection{Equilibrium constant calibration}\label{section:two_particle_equilibrium_constant}

Let $\bar{\pi}_{C} := \bar{\pi}_{(0,0,1)}$ and $\bar{\pi}_{AB} := \bar{\pi}_{(1,1,0)}$ in the notation of Section~\ref{section:continuum_detailed_balance}. In this section we summarize how these probabilities relate to the continuum model equilibrium constant $\Kd$, and discuss how the equilibrium constant $\Kd^h$ for a discretized CRDME model can be chosen to preserve this relationship in the case of a minimal set of substrates for the reversible $A + B \rightleftharpoons C$ reaction, i.e., one $A$ and one $B$ in the unbound state and one $C$ in the bound state. We follow a similar approach for the new CRDME with interaction potentials to that used for calibrating the Brownian dynamics method of~\cite{frohner_reversible_2018}.

In the spatially continuous PBSRD model, with the choice~\eqref{eq:Noe-DB-boltzmann}, we have that 
\begin{equation*}
     \frac{\nu^-(\mathbf{q},\mathbf{q}^+)}{\nu^+(\mathbf{q}^+,\mathbf{q})}  e^{[\Psi^{\mathbf{s}}(\mathbf{q}) - \Psi^{\mathbf{s}^+}(\mathbf{q}^+)]} = 1.
\end{equation*}
The detailed balance relation \eqref{eq:DB-full-equation} then simplifies to 
\begin{equation}\label{eq:multiparticle_DB_relation_partition_function_pis}
\frac{c+1}{a b}\frac{\bar{\pi}_{\mathbf{s}^+}}{\bar{\pi}_{\mathbf{s}}} = \Kd\frac{Z_{\mathbf{s}^+}}{Z_{\mathbf{s}}}.    
\end{equation} 
In the two-particle case, $\mathbf{s} = (a,b,c) = (1,1,0)$, and \eqref{eq:multiparticle_DB_relation_partition_function_pis} reduces to
\begin{equation}\label{eq:DB-twoparticle}
    \frac{\bar{\pi}_{C}}{\bar{\pi}_{AB}} = \Kd \frac{Z_{C}}{Z_{AB}},
\end{equation} where $Z_{C} = Z_{(0,0,1)}$ and $Z_{AB} = Z_{(1,1,0)}$ are the partition functions. From \eqref{eq:DB-twoparticle}, one deduces that to achieve a particular equilibrium ratio $\bar{K} := \frac{\bar{\pi}_{C}}{\bar{\pi}_{AB}}$ of bound to unbound probabilities, one should choose $\Kd = \frac{Z_{AB}}{Z_C}\bar{K}$, as discussed in \cite{frohner_reversible_2018}. %

Our new CRDME with interaction potentials can be shown to be consistent with discrete analogues of the continuum detailed balance relations, including having a unique discrete equilibrium Gibbs-Boltzmann spatial distribution of particles on each stoichiometric compatibility class; see Section~\ref{section:derivations_reactions}. There we show for the general multiparticle CRDME that the probabilities of being in a given number state at equilibrium will satisfy
\begin{equation} \label{eq:multipart_disc_Keq_relation}
    \frac{c+1}{a b} \frac{\bar{\pi}_h^{\bs^+}}{\bar{\pi}_h^{\bs}} =  \Kd^h \frac{Z_h^{\bs^+}}{Z_h^{\bs}},
\end{equation}
where the $h$ subscripts and superscripts indicate the discrete analogues of the continuum quantities and $\Kd^h$ 
is a yet-to-be-chosen mesh-dependent equilibrium constant. Substitution into~\eqref{eq:multpart_crdme_db_step2}, specialized to the two-particle case $\bs = (1,1,0)$ (for which $\alpha^\pm = 1$ and $\hat{\kappa}^\pm = \tilde{\kappa}^\pm$), gives the simplified detailed balance relation that
\begin{equation}\label{eq:discrete_two_particle_DB}
    |V_{ij}|e^{-\tilde{\psi}^{AB}(x_i,x_j)}\tilde{\kappa}^+_{ijk} = \Kd^h|V_{k}|e^{-\psi^{C}(x_k)}\tilde{\kappa}^-_{kij}.
\end{equation} 

To determine a value for $\Kd^h$ in this two-particle setting, we let $\bar{\pi}^h_C$ and $\bar{\pi}^h_{AB}$ be the analogues of $\bar{\pi}_C$ and $\bar{\pi}_{AB}$ in the discrete system and 
\begin{align}\label{eq:discrete_two_particle_partition_functions}
    Z_C^h &:= Z^{(0,0,1)}_h = \sum_{i=1}^N e^{-\psi^C(x_i)}|V_i|, & 
    Z_{AB}^h &:= Z^{(1,1,0)}_h = \sum_{i,j=1}^N e^{-\tilde{\psi}^{AB}(x_i,x_j)}|V_{ij}|
\end{align}
the analogues of the continuum partition functions $Z_C$ and $Z_{AB}$. Then,
\eqref{eq:multipart_disc_Keq_relation} gives that
\begin{align}\label{eq:two_particle_averaged_rates}
    \Kd^h \frac{Z_C^h}{Z_{AB}^h} %
    = \frac{\bar{\pi}_C^h}{\bar{\pi}_{AB}^h}.
\end{align}
It follows from~\eqref{eq:two_particle_averaged_rates} that setting $\Kd^h = \frac{Z^h_{AB}}{Z^h_C}\bar{K}$ results in $\bar{K} = \frac{\bar{\pi}_{C}^h}{\bar{\pi}_{AB}^h}$, %
thereby giving consistency with the continuum PBSRD model's equilibrium constant. Note also that as $h \to 0$, $\Kd^h \to \Kd$ since $Z^h_{AB} \to Z_{AB}$ and $Z^h_C \to Z_C$.

Note that since in practice we only consider short-range interaction potentials, $Z^h_{AB}$ can be computed efficiently as $Z^h_{A} Z^h_{B}$ plus a correction term that only requires summing over pairs of voxels within the interaction distance, where $Z^h_{A} := Z^{(1,0,0)}_h$ and $Z^h_{B} := Z^{(0,1,0)}_h$ are defined analogously to $Z_C^h$ in \eqref{eq:discrete_two_particle_partition_functions}.

\begin{remark}
In the purely diffusive case (see, e.g., \cite{isaacson_unstructured_2018}) it suffices to choose $\Kd^h = \bar{K} |\Omega|$, which is consistent with taking $\Kd = \bar{K} |\Omega|$ in the underlying PBSRD model~\cite{zhang_detailed_2022}. In the presence of one-body potentials, the authors of~\cite{isaacson_zhang_crddme_2024} took a different approach, first computing $\tilde{\kappa}^+_{ijk}$ using a modified version of the formula \eqref{eq:DB-rates-discrete-two-particle} and then directly computing $\tilde{\kappa}^-_{kij}$ from \eqref{eq:discrete_two_particle_DB} to recover the $\Kd$ of the underlying PBSRD model. 
\end{remark}

\subsection{Summary of new CRDME rates}\label{section:summary}

Although the rates given in Section~\ref{section:simulation_algorithm_pair_potentials_reactions} and derived in Section~\ref{section:derivations_reactions_multiparticle} only directly apply to reversible binding/unbinding reactions, the rates for other types of reactions can be derived using the same techniques as in Section~\ref{section:derivations_reactions_multiparticle}. In Table~\ref{tab:discrete_reaction_rates}, we give rates for representative zeroth-, first-, and second-order reactions in our new CRDME. In the first column, $A$, $B$, and $C$ represent different particle types and the subscripts denote the voxel location; for example, $A_i\to B_j$ denotes the conversion of an $A$ in voxel $V_i$ to a $B$ in voxel $V_j$. The base rates in the second column represent rates for irreversible reactions. In the third column, we give \textit{modulating} (two-particle Boltzmann) factors in the presence of potentials. For \textit{dilute} systems that contain exactly as many substrates as needed for one reaction instance to be possible, e.g., one $A$ and $B$ pair or one $C$ for $A + B \rightleftharpoons C$, multiplying these factors into the base rates yields rates satisfying detailed balance. Finally, for general multiparticle reactions we give acceptance (multiparticle Boltzmann) factors accounting for interactions with particles not directly involved in the reaction. In Table~\ref{tab:discrete_reaction_rates}  we use $\mathbf{q}$ to represent a generic pre-event state and $\mathbf{q}^R$ to represent the post-event state. 

\begin{table}[tb]
    \centering
     \resizebox{\textwidth}{!}{
    \begin{tabular}{l|l|l|l}
       Event type & Base rate & Modulating factor & Acceptance factor \\\hline
        $\varnothing \to A_i$ & $\kappa_i |V_i|$ & $e^{-\psi^A(x_i)}$ & $\min\{1,e^{-[\Psi(\mathbf{q}^R) - \Psi(\mathbf{q})] + \psi^A(x_i)}\}$ \\ 
        $A_i \to \varnothing$ & $\kappa_i a_i$ & --- & $\min\{1,e^{-[\Psi(\mathbf{q}^R) - \Psi(\mathbf{q})] - \psi^A(x_i)}\}$ \\ 
        $A_i\to B_j$ & $\kappa_i\delta_{ij} a_i$ & $e^{-\psi^B(x_j)}$ & $\min\{1,e^{-[\Psi(\mathbf{q}^R) - \Psi(\mathbf{q})] + \psi^B(x_j) - \psi^A(x_i)}\}$\\
        $A_i + B_j \to C_k$ &  $\kappa_{ijk}^+ a_i b_j$ & $e^{-\psi^C(x_k)}$ & $\min\{1,e^{-[\Psi(\mathbf{q}^R) - \Psi(\mathbf{q})] + \psi^C(x_k) - \tilde{\psi}^{AB}(x_i,x_j)}\}$\\
        $C_k \to A_i + B_j$ & $\kappa^-_{kij} c_k$ & $e^{-\tilde{\psi}^{AB}(x_i,x_j)}$ & $\min\{1,e^{-[\Psi(\mathbf{q}^R) - \Psi(\mathbf{q})] + \tilde{\psi}^{AB}(x_i,x_j) - \psi^C(x_k)}\}$\\
        $A_i \to A_j$ & $D^A\frac{\omega_{ji}}{|V_i|} a_i$ & $B(\Psi(\mathbf{q}^R) - \Psi(\mathbf{q}))$ & ---
    \end{tabular}
    }
    \caption{CRDME hopping and reaction rates for irreversible and reversible reactions. $a_i$ denotes the number of particles of species $A$ in voxel $V_i$, with $b_j$ and $c_k$ defined analogously. For transitions representing chemical reactions, the base rate is the reaction rate for irreversible reactions, while the product of the base rate and the modulating factor gives the proposal rate for reversible reactions. For the latter case, the acceptance factor then gives the probability of accepting a proposed reaction. Note that the factors given for the $\varnothing \rightleftharpoons A$ reactions are only needed when $\varnothing$ is interpreted as an equilibrium particle reservoir with which the system exchanges particles. The last row represents spatial hopping transitions, for which the base rate is the rate in the pure diffusion case, while the modulating factor multiplies it to provide a rate that accounts for potential interactions. In rows 1--3, $\kappa_i$ is the value in voxel $V_i$ of a chosen microscopic rate function $\kappa:\Omega\to [0,\infty)$, which has units of (volume\,$\cdot$\,time)$^{-1}$ for the zeroth-order reaction and time$^{-1}$ for the first-order reactions. For smooth $\kappa$ we take $\kappa_i = \kappa(x_i)$, while for less smooth rates (e.g., $\kappa \in L^\infty(\Omega)$) we could instead choose the volume average $\kappa_i = \frac{1}{|V_i|} \int_{V_i} \kappa(x)\,dx$. In rows 4--5, $\kappa^+_{ijk}$ and $\kappa^-_{kij}$ refer to any of \eqref{eq:discrete_fwd_standard}--\eqref{eq:discrete_bwd_smooth}.}
\label{tab:discrete_reaction_rates}
\end{table}

As we mentioned at the end of Section~\ref{section:simulation_algorithm_pair_potentials_reactions}, to apply the rates from Table~\ref{tab:discrete_reaction_rates} in a CRDME simulation there are two cases for reaction events:

\begin{enumerate}
    \item Irreversible reactions occur at the base rate.
    \item Reversible reactions are proposed at the \textit{proposal rate} (product of the base rate with the modulating factor). Proposed reactions are \textit{accepted} with probability given by the acceptance factor.
\end{enumerate}

%% file: manuscript/results.tex
\section{Numerical results}\label{section:numerics_pair_potentials}

In this section, we present validation and convergence tests for the CRDME simulation method presented earlier by computing statistics for simulations of different particle systems. We start by computing first reaction times in the $A+B\to\varnothing$ reaction with Lennard-Jones interaction potentials. Then, we measure bound state probabilities for a two-particle reversible reaction $A + B \rightleftharpoons C$ with harmonic potentials between all pairs of particles. Finally, we validate our method by studying a multiparticle version of the same reversible system, comparing our results with the Brownian dynamics-based ReaDDy software \cite{hoffmann_readdy_2019}. %

\subsection{Annihilation reaction}\label{section:results_annihilation}

The first system we study involves two particle types $A$ and $B$, which can undergo the annihilation reaction $$A + B \xrightarrow[]{\,\kappa(x,y)\,} \varnothing,$$ via the Doi reaction kernel \eqref{eq:doi_kernel}. The interaction distance $\varepsilon$ can be viewed as the sum of the particle radii. To compute the discrete reaction rate $\kappa_{ij}$ for an $A_i + B_j \to \varnothing$ reaction, we take $\gamma=1$ in \eqref{eq:discrete_fwd_lump} and set $\kappa_{ij} = \kappa^+_{iji}$, i.e., 
\begin{equation*}
  \kappa_{ij} = \frac{\lambda}{|V_j|} \int_{V_j} 1_{|x_i - y| < \varepsilon}dy. 
\end{equation*}

In addition to the annihilation reaction, the particles interact via a singular 6-12 Lennard-Jones potential. The Lennard-Jones potential is repulsive at distances below $2^{1/6}\varepsilon$ (in particular when the particles contact each other, i.e., at distances below $\varepsilon$), weakly attractive when $|x-y|$ is within the range $(2^{1/6}\varepsilon,\varepsilon_{\text{cutoff}})$, and zero beyond $\varepsilon_{\text{cutoff}}$. Writing $r = |x-y|$ for the distance between the two particles, the formula is given by
$$
\psi_{\text{LJ}}(r) =
\begin{cases} 
\eta_{S_1 S_2} \left(\left(\frac{\varepsilon}{r}\right)^{12} - \left(\frac{\varepsilon}{r}\right)^6\right) & r \leq \varepsilon_{\text{cutoff}},\\
0 & r > \varepsilon_{\text{cutoff}}.
\end{cases}
$$
The cutoff radius for the potential is usually chosen as a multiple of the particle radii. Here we use ${\varepsilon_{\text{cutoff}} = \frac{5}{2}\varepsilon}$.

We also introduce two modifications to $\psi_{\text{LJ}}$ as presented above. The first is a shift parameter $\psi_{\text{shift}} := \psi_{\text{LJ}}(\varepsilon_{\text{cutoff}})$, which modifies $\psi_{\text{LJ}}(r)$ to $\psi_{\text{LJ}}(r) - \psi_{\text{shift}}$, ensuring that $\psi_{\text{LJ}}$ is continuous. The second modification concerns the singularity, which is important in our context since our simulation method involves evaluation of the potential at nodes of the mesh. Since $\psi_{\text{LJ}}(0) = \infty$, one way to deal with the singularity would be to directly disallow placement of two particles with a Lennard-Jones potential into the same voxel and set the hopping rates into occupied voxels to zero. In our simulations, we found it simpler to modify $\psi_{\text{LJ}}$ to a mesh-dependent $\psi_{\text{LJ}}^{h}$ with $\psi_{\text{LJ}}^{h}(0) = \psi_{\text{LJ}}\left(\underline{r}\right) - \psi_{\text{shift}},$ where we choose $\underline{r}$ to be a small fraction of the minimum mesh spacing $h_{\text{min}}$, ensuring that the hopping rates to occupied voxels are effectively zero. After experimentation we chose the final value $\underline{r} = \frac{h_{\text{min}}}{8}$. Our final formula is thus
$$
\psi_{\text{LJ}}^h(r) =
\begin{cases}
\psi_{\text{LJ}}(\underline{r}) - \psi_{\text{shift}} & r < \underline{r},\\
\eta_{S_1 S_2} \left(\left(\frac{\varepsilon}{r}\right)^{12} - \left(\frac{\varepsilon}{r}\right)^6\right) - \psi_{\text{shift}} &\underline{r}\leq r < \varepsilon_{\text{cutoff}},\\
0 & r \geq \varepsilon_{\text{cutoff}}.
\end{cases}
$$
The full list of parameters we use in our simulations is shown in Table~\ref{tab:LJ-sim-parameters}. 

\begin{table}[tb]
    \centering
    \begin{tabular}{c|c|c}
        \hline
        Parameter & Description & Value \\\hline
        $D^A$ & Diffusion constant, particle type $A$ & $10^3$ \\\hline
        $D^B$ & Diffusion constant, particle type $B$ & $10^3$ \\\hline
        $\rho_A$ & Initial distribution, particle type $A$ & $e^{-\frac{(x - .5)^2 + (y-.5)^2}{.02}}$ \\\hline
        $\rho_B$ & Initial distribution, particle type $B$ & $e^{-\frac{(x + .5)^2 + (y + .5)^2}{.02}}$ \\\hline
        $r_A$ & Particle $A$ radius & .02 \\ \hline
        $r_B$ & Particle $B$ radius & .02 \\\hline
        $\varepsilon$ & Doi reaction radius/Lennard-Jones potential parameter & .04 \\\hline
        $\lambda$ & Doi reaction rate & $10^6$ \\\hline
        $\varepsilon_{\text{cutoff}}$ & Lennard-Jones potential cutoff & $\frac{5}{2}\varepsilon = .1$ \\\hline
        $\eta_{S_1 S_2}$ & Lennard-Jones potential interaction strength & 1 \\\hline
        $\underline{r}$ & Lennard-Jones singularity cutoff distance & $\frac{h_{\text{min}}}{8}$ \\\hline
    \end{tabular}
    \caption{Parameters used in $A+B\to\varnothing$ simulations with a Lennard-Jones potential.
    }
    \label{tab:LJ-sim-parameters}
\end{table}

We start each simulation with $A_0$ particles each of types $A$ and $B$, letting $A^{h,\ell}(t)$ be the number of $A$ particles remaining at time $t \geq 0$. We initially distribute the particles using normalized versions of the sharp Gaussian distributions $\rho_A$ and $\rho_B$ listed in Table~\ref{tab:LJ-sim-parameters}, which we chose to have minimal overlap to avoid initializing $A$ and $B$ particles in the same voxels.

We measure the proportion of particles $\frac{A^{h,\ell}(t)}{A_0}$ remaining at a set of times $t = t_0=0,t_1,\dots,t_K=0.1$. We also record the time taken for all particles to react for each trial 
\begin{equation*}
  T_{\text{react}}^{h,\ell} = \inf \{t \geq 0 : A^{h,\ell}(t) = 0\},
\end{equation*}
and measure convergence of the mean value of this statistic over $P$ trials, 
$$T_{\text{react}}^{h} = \frac{1}{P}\sum_{\ell=1}^P T_{\text{react}}^{h,\ell}.$$ 
For the two-particle models, we also plot the discrete survival probability $P(T_{\text{react}}^{h,\ell} > t)$. Since we do not have an analytical solution for this problem, we measure the convergence rate of $T_{\text{react}}^{h}$ empirically by looking at differences $$|T_{\text{react}}^{h} - T_{\text{react}}^H|,\qquad H\approx 2h,$$ over consecutive mesh refinements.

\subsubsection{Two-particle case}\label{section:results_annihilation_two_particle}

We simulate the two-particle problem on two different two-dimensional domains, $\Omega_{\text{square}} = [-1,1]^2$ and $\Omega_{\text{disk}} = B_1(0) = \{x \in \RR^2 : |x| < 1\}$. On $\Omega_{\text{square}}$, we generate a hierarchy of meshes by starting with an initial mesh consisting of only two triangles and recursively subdividing each triangle from the previous mesh into four triangles by connecting the midpoints of the coarser mesh edges. For $\Omega_{\text{disk}}$, we construct our hierarchy of meshes by retriangulating the domain with different target element sizes, approximating the circular boundary by polygons with increasing numbers of vertices. Approximating the boundary of the disk by polygons introduces additional error on top of the usual FEM approximation error; however, this error does not affect our empirical convergence rates. 

\begin{figure}[ht]
    \centering
\includegraphics[scale=.4]{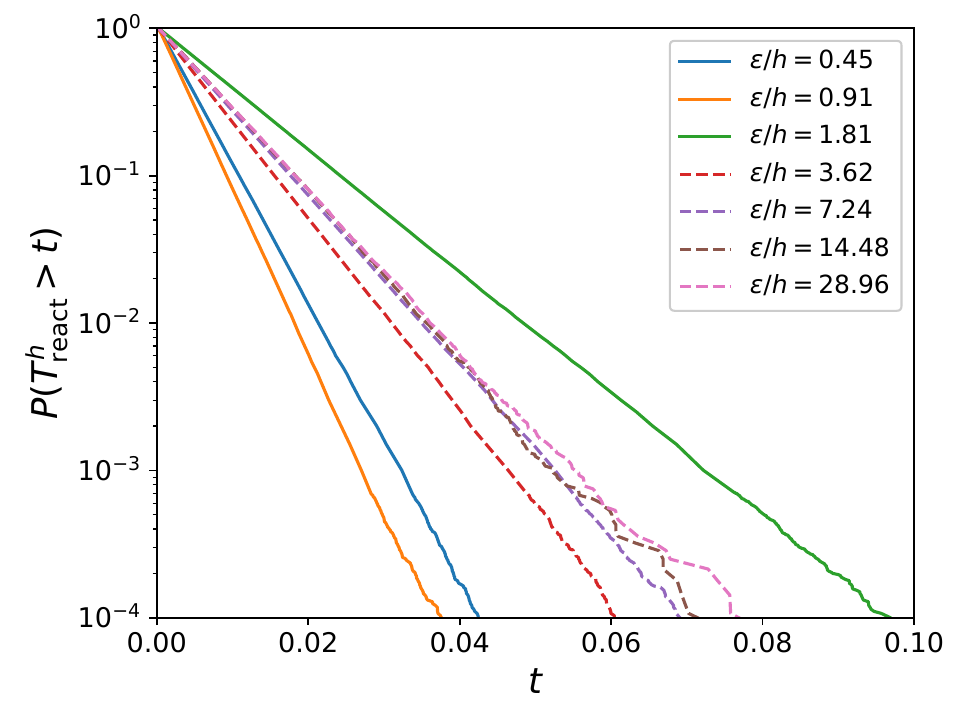}\includegraphics[scale=.4]{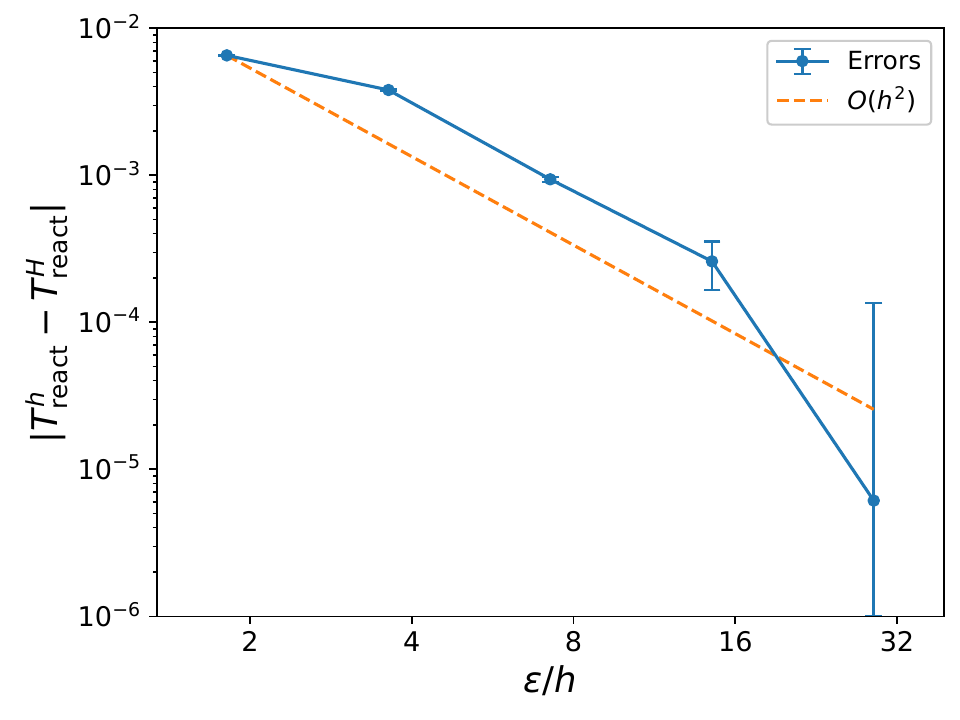}
\includegraphics[scale=.4]{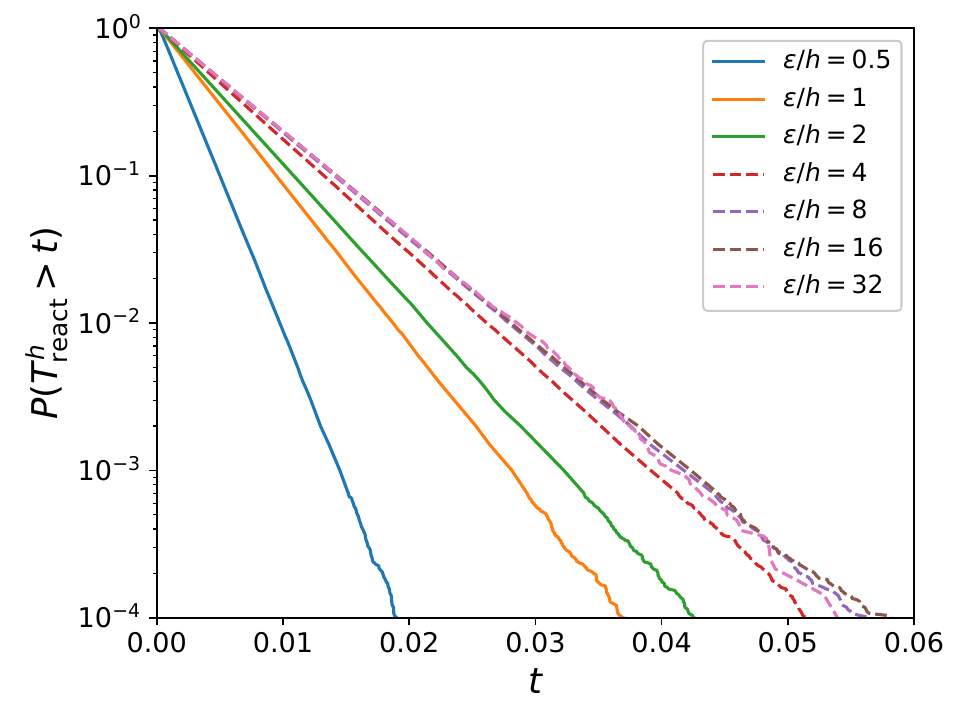}\includegraphics[scale=.4]{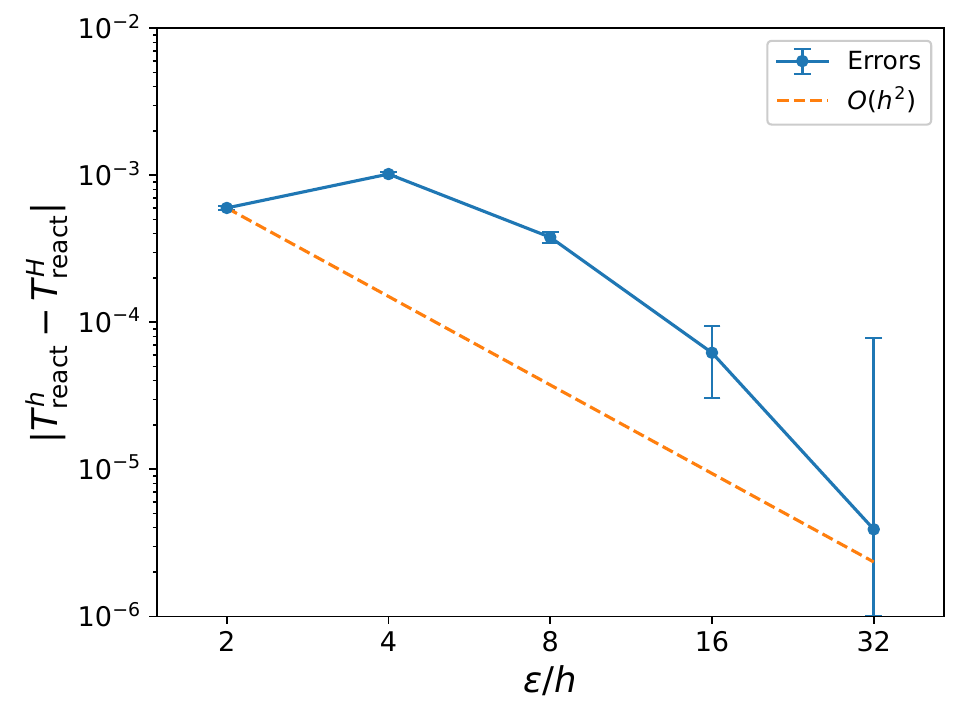}
    \caption{Discrete survival probabilities for different mesh sizes (left column) and consecutive differences in mean reaction time $|T_{\text{react}}^h - T_{\text{react}}^H|$, where $H \approx 2h$ (right column), for the two-particle Lennard-Jones model with an annihilation reaction on a square domain (top row) and a disk domain (bottom row). Error bars in the right column show 95\% confidence intervals for the differences, computed from the standard errors of the two means; a bar that reaches the bottom of the plot indicates an interval that includes zero.}
    \label{fig:two-particle-annihilation}
\end{figure}

In Figure~\ref{fig:two-particle-annihilation}, we show the discrete survival probabilities for the first reaction time (left column) and the convergence of the mean first reaction times (right column) on $\Omega_{\text{square}}$ (top row) and $\Omega_{\text{disk}}$ (bottom row). The values are averaged over $P = 280\text{,000}$ simulations, except on the finest two square meshes and the finest disk mesh, where $P = 28\text{,000}$. For the square domain, the differences are consistent with $O(h^2)$ convergence within the 95\% confidence intervals once $\frac{\varepsilon}{h} \gtrsim 4$; the error bar on the last difference is wide, but its 95\% confidence interval contains the value obtained by extrapolating the previous two differences at rate $O(h^2)$. On the disk, the differences are likewise consistent with $O(h^2)$ convergence once $\frac{\varepsilon}{h} \gtrsim 8$, with the last confidence interval again containing the $O(h^2)$ extrapolation of the previous two differences. The survival probabilities for the reaction times also exhibit good convergence, indicating that we capture the full distribution of reaction times (and not just the means) correctly. On our two finest mesh refinements for both the disk and square meshes, the survival probabilities are essentially overlapping.

\begin{figure}[h]
    \centering
\includegraphics[scale=.4]{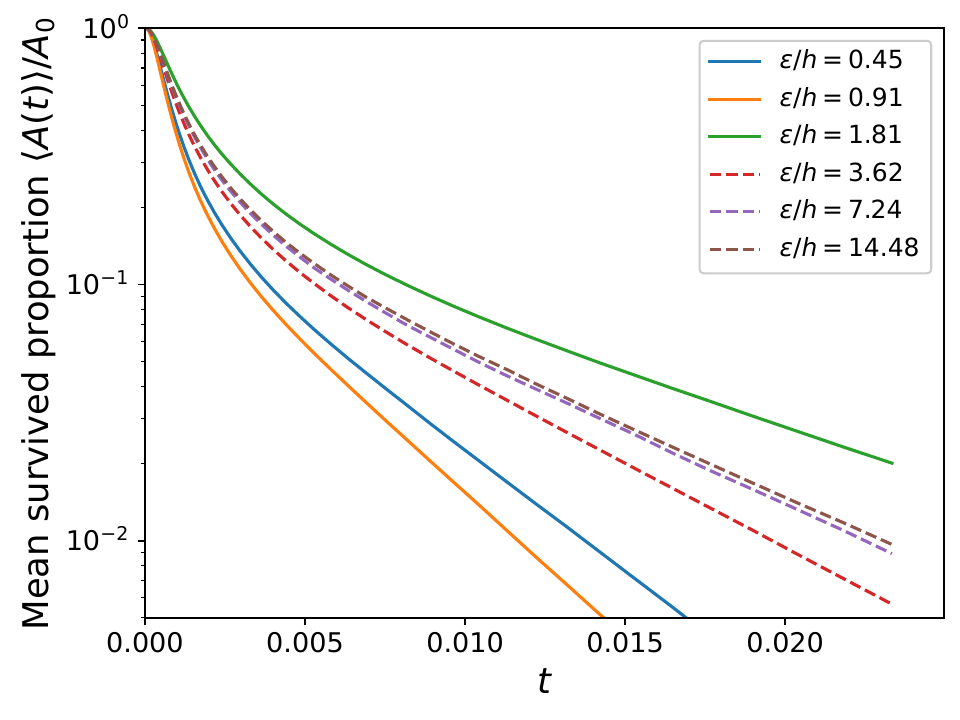}\includegraphics[scale=.4]{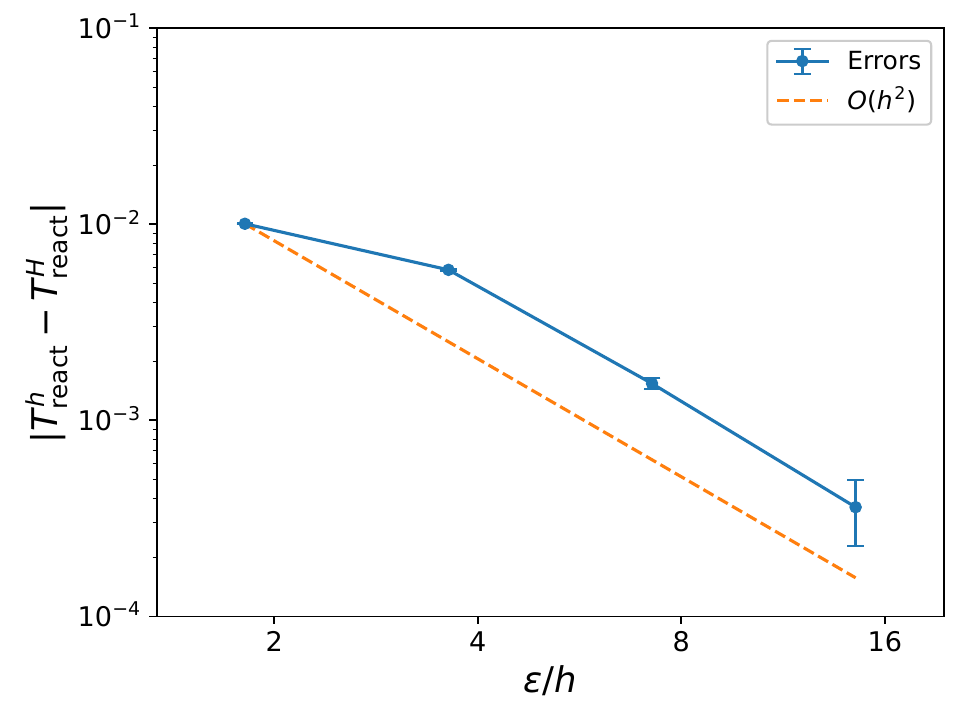}
\includegraphics[scale=.4]{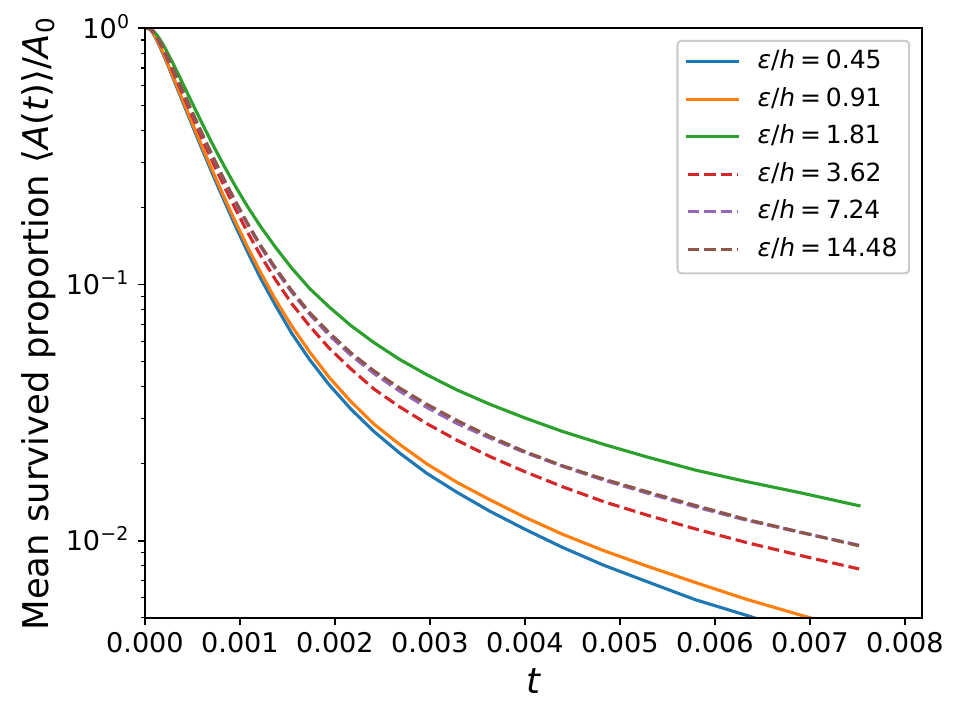}\includegraphics[scale=.4]{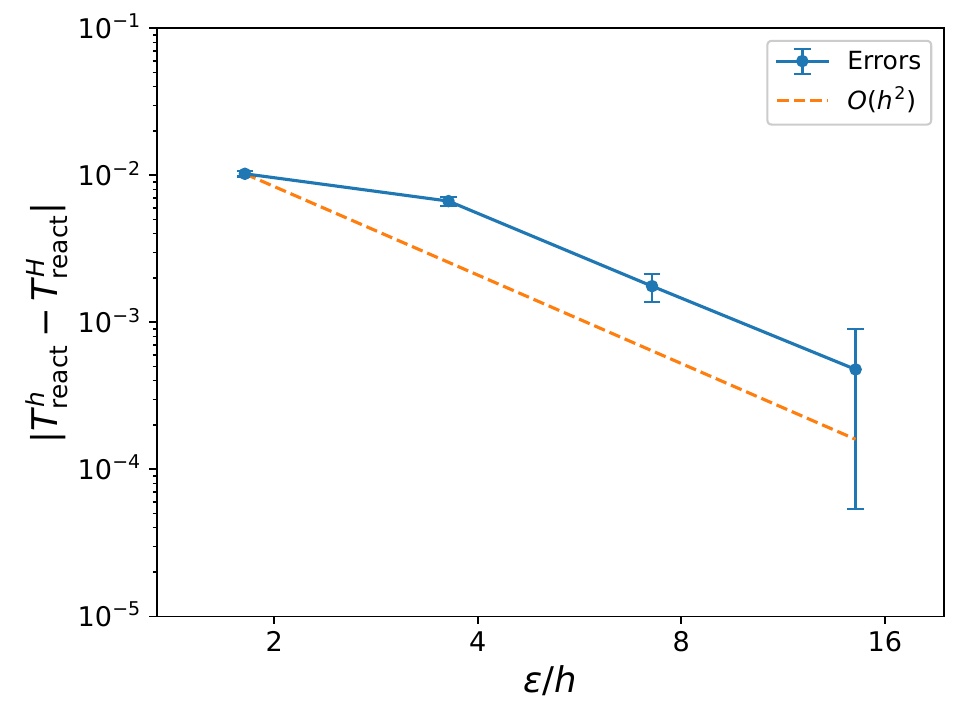}
    \caption{Proportion of particles remaining for different mesh sizes (left column) and consecutive differences in mean reaction time $|T_{\text{react}}^h - T_{\text{react}}^H|$, where $H \approx 2h$ (right column), for the multiparticle Lennard-Jones model with an annihilation reaction on a square domain with starting population $A_0 = 10$ (top row) and $A_0 = 100$ (bottom row). Error bars show 95\% confidence intervals for the differences, as in Figure~\ref{fig:two-particle-annihilation}.}
    \label{fig:multiparticle-annihilation}
\end{figure}
\subsubsection{Multiparticle case}\label{section:results_annihiliation_multiparticle}

We simulate the multiparticle case on $\Omega_{\text{square}}$ only, plotting in Figure~\ref{fig:multiparticle-annihilation} the mean proportion of particles remaining over time for $A_0 = B_0 \in \{10,100\}$ (left column) and the difference in mean reaction time on consecutive meshes (right column). The mean values are averaged over $P=280\text{,000}$ trials for the 10-particle model except on the finest two meshes, where we use $P=28\text{,000}$. For the 100-particle model we use $P = 2\text{,800}$ trials on every mesh. We use fewer trials on the finest meshes because of the high computational expense --- the work required (in two dimensions) scales roughly like $O(\varepsilon^2 Nh^{-2}) = O(\varepsilon^2 h^{-4})$, since the number of diffusion events is $O(h^{-2})$ and in our implementation, each $B$-particle diffusion event requires $O(\varepsilon^2 N)$ reaction rate updates.

From the results in Figure~\ref{fig:multiparticle-annihilation}, we again see convergence rates consistent with $O(h^2)$ within the 95\% confidence intervals for the mean reaction times once $\frac{\varepsilon}{h} \gtrsim 4$, and good convergence for the mean proportion of particles surviving over time as well.

\subsection{Reversible reaction}\label{section:results_reversible}

In this section, we study the full version of our model problem by adding a reversible reaction. We start with convergence tests for a two-particle model on the disk meshes of Section~\ref{section:results_annihilation_two_particle}, rescaled as described below. We then compare our results on a three-dimensional structured-grid multiparticle model on a periodic cubic domain with fully reversible reaction dynamics to the Brownian dynamics-based ReaDDy software \cite{hoffmann_readdy_2019}.

For both models, we use harmonic interaction potentials between each pair of particles of the form 
\begin{equation}\label{eq:harmonic_potential}
    \psi^{S_1S_2}(r) = \frac{\eta_{S_1 S_2}}{2}\max\{(r_{S_1} + r_{S_2}) - r,0\}^2,
\end{equation} 
where $S_1,S_2\in\{A,B,C\}$ and $r_{S_1}$ and $r_{S_2}$ are the radii of the particles, so that the potential acts only within the interaction radius $\varepsilon_{S_1S_2} := r_{S_1} + r_{S_2}$. We set $\eta_{S_1 S_2} = 0.1$ in the two-particle simulations and $\eta_{S_1 S_2} = 10$ in the multiparticle simulations. For the reaction rates, we again use the Doi reaction kernel \eqref{eq:doi_kernel} within the discrete reaction base rates \eqref{eq:discrete_fwd_lump} and \eqref{eq:discrete_bwd_lump} with $\gamma = \frac{1}{2}$. For the Doi kernel reaction radius, we choose $\varepsilon = \varepsilon_{AB} = r_A + r_B$. For the periodic cubic domain simulations, the backward rate $\mu$ in Table~\ref{tab:multiparticle_reversible_parameters} is the total dissociation rate of a $C$ particle. Periodic translational invariance, together with the fact that both the potential's range and the reactive interaction distance are small relative to the domain length, gives $\Kd = \frac{\lambda}{\mu}\int_{|r|<\varepsilon}e^{-\psi^{AB}(|r|)}\,dr$. For the disk domain simulations with parameters in Table~\ref{tab:two_particle_reversible_parameters}, the backward rate $\mu$ is instead the total dissociation rate of a $C$ particle in the absence of the pair potential, which corresponds to $\Kd = \frac{\lambda}{\mu}|B_\varepsilon(0)| = \frac{\lambda\pi\varepsilon^2}{\mu}$, with $B_\varepsilon(x)$ denoting the $d$-dimensional ball of radius $\varepsilon$ centered at $x$. In that case, for a $C$ particle at $z$, the reactive region is strictly $B_\varepsilon(z)\cap\Omega$; because both the potential range and the reactive interaction distance are small relative to the disk radius, we expect that replacing this position-dependent intersection by the full ball introduces only a negligible correction; see~\cite{zhang_detailed_2022} for related discussion in the purely diffusive case. In neither geometry do we explicitly calibrate a mesh-dependent $\Kd^h$ as described in Section~\ref{section:two_particle_equilibrium_constant}; we instead discretize the reaction kernels directly from the continuum description. In both geometries we initially distribute the particles using a uniform distribution.

 Let $A_0$, $B_0$, and $C_0$ be the number of type $A$, $B$, and $C$ particles initially present in the system. We always choose $A_0=B_0$ in our simulations, so that $M_0 = A_0 + C_0$ is the maximum number of $C$ particles, or bound $A$-$B$ particle pairs, present in a given simulation. Let $\bs_c = (M_0-c,M_0-c,c)$ denote the number state with $c$ particles of type $C$. We estimate the mean total proportion of particles in the bound state 
 \begin{equation*}
 P_\text{bound}(t) = \avg{\frac{C(t)}{M_0}} = \sum_{c=0}^{M_0} \frac{c}{M_0 ((M_0-c)!)^2 c!}\int_{\Omega^{2M_0-c}} f^{\bs_c}(\mathbf{q},t)\, d\mathbf{q}
 \end{equation*}
 using the estimator $$P^{h}_{\text{bound}}(t) = \frac{1}{P}\sum_{\ell=1}^P \frac{C^{h,\ell}(t)}{M_0},$$ where $C^{h,\ell}(t)$ is the total number of particles of type $C$ in the system at time $t$ during trial $\ell$.

 \subsubsection{Two-particle case}\label{section:results_reversible_two_particle}

For the two-particle simulation, we use the disk domain and associated meshes from Section~\ref{section:results_annihilation_two_particle}, rescaled so that the disk has radius $R = 297.363$, with the dimensionless parameters listed in Table~\ref{tab:two_particle_reversible_parameters}. Each trial starts with a single $C$ particle ($A_0 = B_0 = 0$, $C_0 = M_0 = 1$) placed uniformly in the disk.

\begin{table}[tb]
    \centering
\begin{tabular}{c|c|c}
\hline
   Parameter & Description & Value \\ \hline
   $D^A$  & Diffusion constant, species $A$ & 2862  \\\hline  
     $D^B$  & Diffusion constant, species $B$ & 1522 \\  \hline
$D^C$  & Diffusion constant, species $C$ & 1336 \\ \hline 
$r_A$ & Particle $A$ radius & 1.875 \\\hline
$r_B$ & Particle $B$ radius & 3.75 \\\hline
$r_C$ & Particle $C$ radius & 7.8 \\\hline
$R$ & Disk radius & 297.363 \\\hline
$\varepsilon_{S_1S_2}$ & Potential interaction radii & $r_{S_1} + r_{S_2}$ \\\hline
$\eta_{S_1S_2}$ & Harmonic potential interaction strength & 0.1 \\\hline
$\lambda$ & Microscopic Doi forward rate & $10^2$ \\\hline
$\mu$ & Total dissociation rate of a $C$ particle without the pair potential & 28.454 \\\hline
$\gamma$ & Reaction product placement probability & $\frac{1}{2}$ \\\hline

$\rho_S(x)$ & Initial condition, species $S$ & Uniform \\\hline
\end{tabular}
    \caption{Parameters for the two-particle reversible model.}
    \label{tab:two_particle_reversible_parameters}
\end{table}

\begin{figure}[h]
    \centering
            \includegraphics[scale=.4]{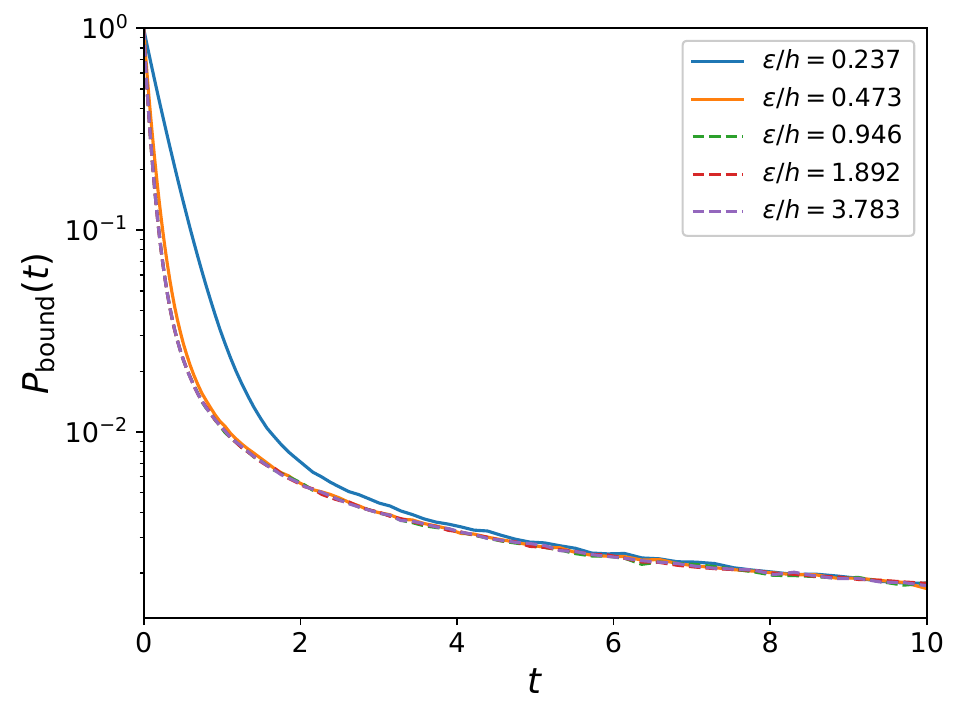}\includegraphics[scale=.4]{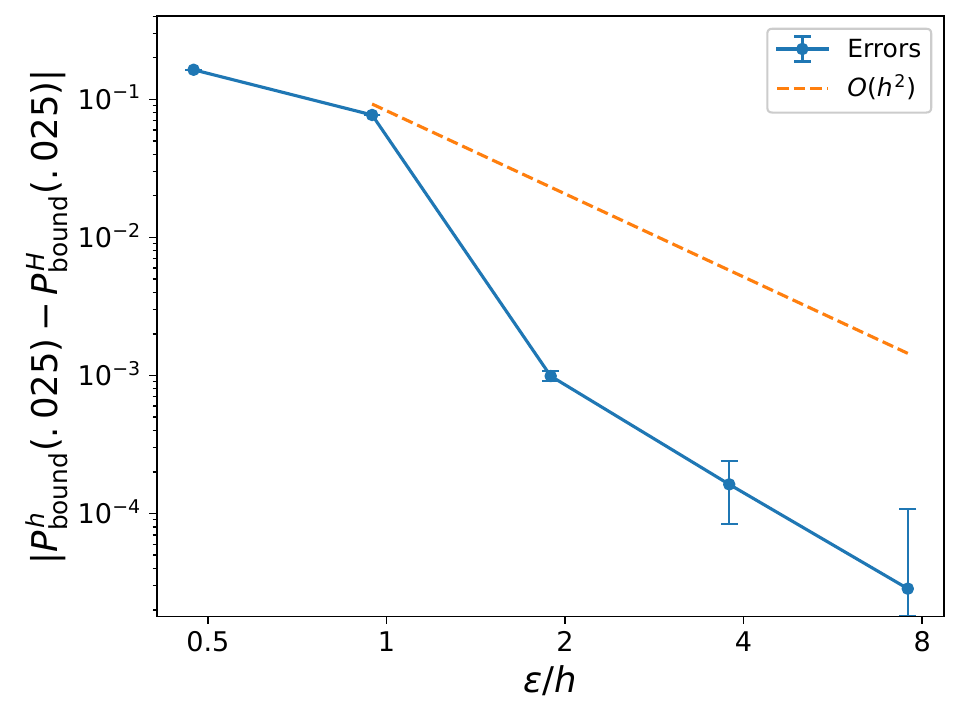}
    \caption{Empirical estimate of $P_{\text{bound}}(t)$ over time (left) and consecutive differences $|P_{\text{bound}}^h(.025) - P_{\text{bound}}^H(.025)|$, where $H \approx 2h$, on different meshes (right), with 95\% confidence intervals for the differences; the bar on the last point reaches the bottom of the plot because its interval includes zero.}
    \label{fig:twoparticle-reversible}
\end{figure}

As in Section~\ref{section:results_annihilation}, we measure the convergence rate of $P_{\text{bound}}^{h}(t)$ at the fixed time $t = .025$ empirically by looking at differences $$|P_{\text{bound}}^{h}(.025) - P_{\text{bound}}^{H}(.025)|,\qquad H\approx 2h,$$ over consecutive mesh refinements. The curves in the left panel of Figure~\ref{fig:twoparticle-reversible} are averaged over $P = 2.8\cdot 10^{6}$ trials per mesh, and the differences in the right panel are computed from separate runs with $P = 2.8\cdot 10^{8}$ trials per mesh. The results in Figure~\ref{fig:twoparticle-reversible} show that once the mesh is refined below the interaction radius $\varepsilon = 5.625$ (i.e., $\frac{\varepsilon}{h} > 1$), the differences drop sharply and then decrease at least as fast as $O(h^2)$; the error bar on the last difference is wide, but its 95\% confidence interval lies below the $O(h^2)$ reference line and contains the value obtained by extrapolating the previous two differences at rate $O(h^2)$.

We also see that the green, red, and purple curves are almost completely overlapping, indicating that the bound state probabilities, plotted as a function of time, converge quickly once $\frac{\varepsilon}{h} \approx 1$; the convergence plot shows a drop by a factor of about $80$ in the consecutive differences once $\frac{\varepsilon}{h} > 1$. We have omitted the last refinement level on the left plot because there is no visual difference between it and the last two curves.

\subsubsection{Multiparticle case}\label{section:results_reversible_multiparticle}

\begin{table}[tb]
    \centering
\begin{tabular}{c|c|c}
\hline
   Parameter & Description & Value \\ \hline
   $D^A$  & Diffusion constant, species $A$ & 5 \\ \hline 
     $D^B$  & Diffusion constant, species $B$ & 5 \\  \hline
$D^C$  & Diffusion constant, species $C$ & 7 \\  \hline
$r_A$ & Particle $A$ radius & 1 \\\hline
$r_B$ & Particle $B$ radius & 1 \\\hline
$r_C$ & Particle $C$ radius & 1.4 \\\hline
$\varepsilon_{S_1S_2}$ & Potential interaction radii & $r_{S_1} + r_{S_2}$ \\\hline
$\eta_{S_1S_2}$ & Harmonic potential interaction strength & 10 \\\hline
$\lambda$ & Microscopic Doi forward rate & .05 \\\hline
$\mu$ & Total dissociation rate of a $C$ particle & .001 \\\hline
$\gamma$ & Reaction product placement probability & $\frac{1}{2}$ \\\hline
$\rho_S(x)$ & Initial condition, species $S$ & Uniform \\\hline
\end{tabular}
    \caption{Parameters for the multiparticle reversible model.}
    \label{tab:multiparticle_reversible_parameters}
\end{table}

In this section, we study the fully reversible reaction of Section~\ref{section:results_reversible_two_particle} in a multiparticle system. We will compare our results to results obtained using the ReaDDy software \cite{hoffmann_readdy_2019}. Since ReaDDy is only implemented on three-dimensional box domains, we employ a standard tetrahedralization of the periodic box domain $\Omega = [0,20]^3$ along with a Cartesian dual mesh (see the discussion at the beginning of Section~\ref{section:simulation_algorithm_background}). We choose the remaining parameters (except the microscopic reaction rates) to match those used in an example from \cite{frohner_reversible_2018} that the authors use to test their Brownian dynamics-based reversible particle dynamics. The parameters are listed in Table~\ref{tab:multiparticle_reversible_parameters}. We note that since the current ReaDDy implementation places product particles for $A+B\to C$ reactions at the midpoint between the substrates, we made a slight modification to their source code to support the alternative placement mechanism \eqref{eq:nodal_placement}.

\begin{figure}[h]
    \centering
    \includegraphics[width=\textwidth]{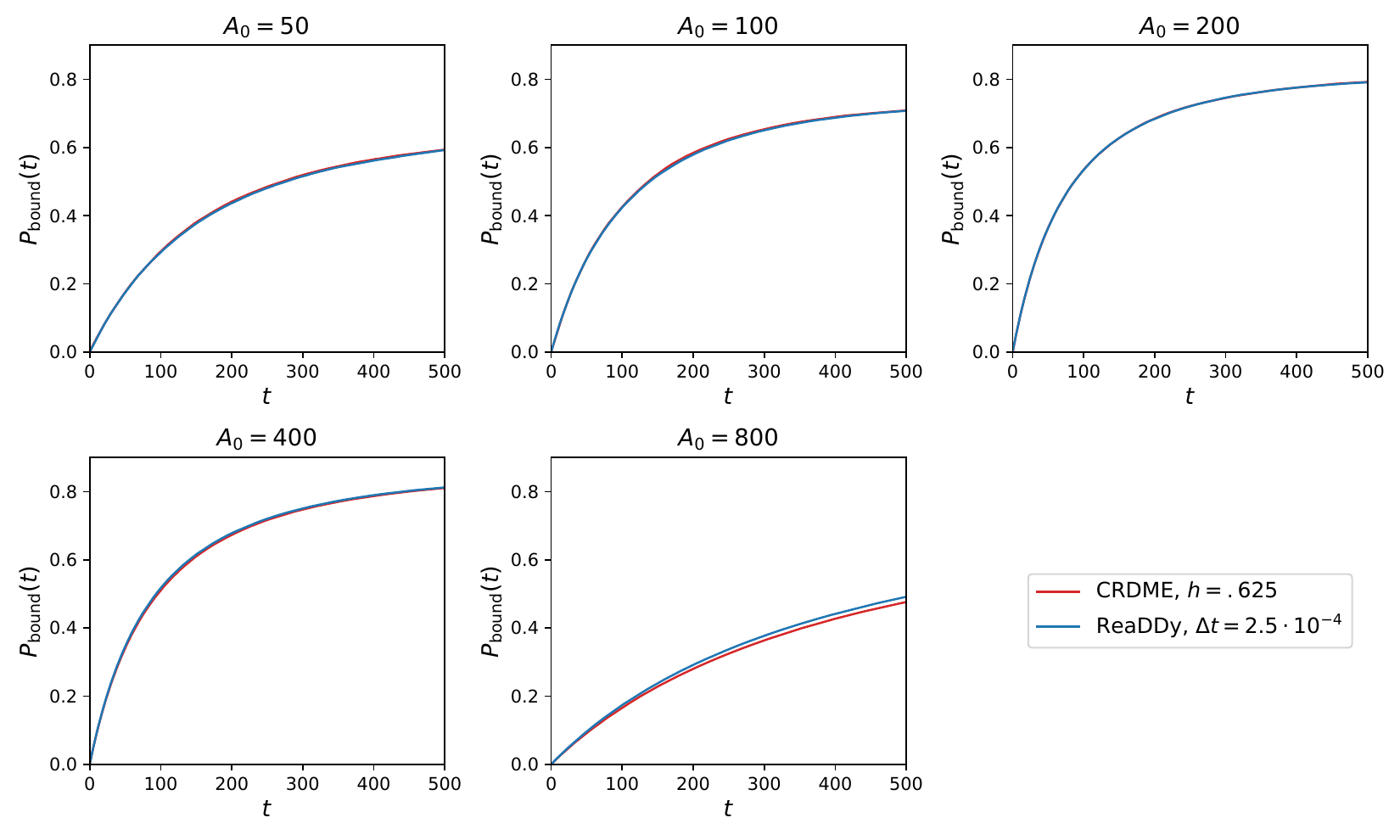}
    
    \caption{Comparison of transients for the proportion of bound particles $P_{\text{bound}}(t) = \frac{\avg{C(t)}}{A_0}$ for ReaDDy and CRDME simulations with $\Delta t_{\text{ReaDDy}} = 2.5\cdot 10^{-4}$ and $h=.625$, respectively.  The results agree closely except in the $A_0 =800$ case. }
    \label{fig:multiparticle_transient_comparision}
\end{figure}

In Figure~\ref{fig:multiparticle_transient_comparision}, we compare the mean values for the transients produced by CRDME and ReaDDy simulations. For both the CRDME and ReaDDy simulations, we started with initial populations of $M_0 = A_0 = B_0 \in \{50,100,200,400,800\}$, $C_0 = 0$, and ran the simulations until $T_f = 500$. For the ReaDDy simulations, we chose the Brownian dynamics time step $\Delta t = 2.5\cdot 10^{-4}$ by taking refinements $2^{-k}\cdot .001$ for $k=0,1,2,\dots$ until the difference between successive values of $\avg{C(t)}$ was less than 1\%, where $\avg{C(t)}$ is the average number of particles of type $C$ at time $t$ in the $A_0 = 50$ simulation. We note that our time step is half the value of $5 \cdot 10^{-4}$ used in \cite{hoffmann_readdy_2019}. Similarly, we chose our CRDME refinement level by taking meshes with $2^L$ (unique) nodes in each direction, with $L$ chosen under the constraint that $h = \frac{20}{2^{L}}$ was smaller than the minimum interaction radius $\varepsilon = 2$ for our model. We found that $L=5$, so that $h=.625 < 2$ and $N=32{,}768$, provided an accurate solution. %

\begin{figure}[h]
    \centering
    \includegraphics[scale=.35]{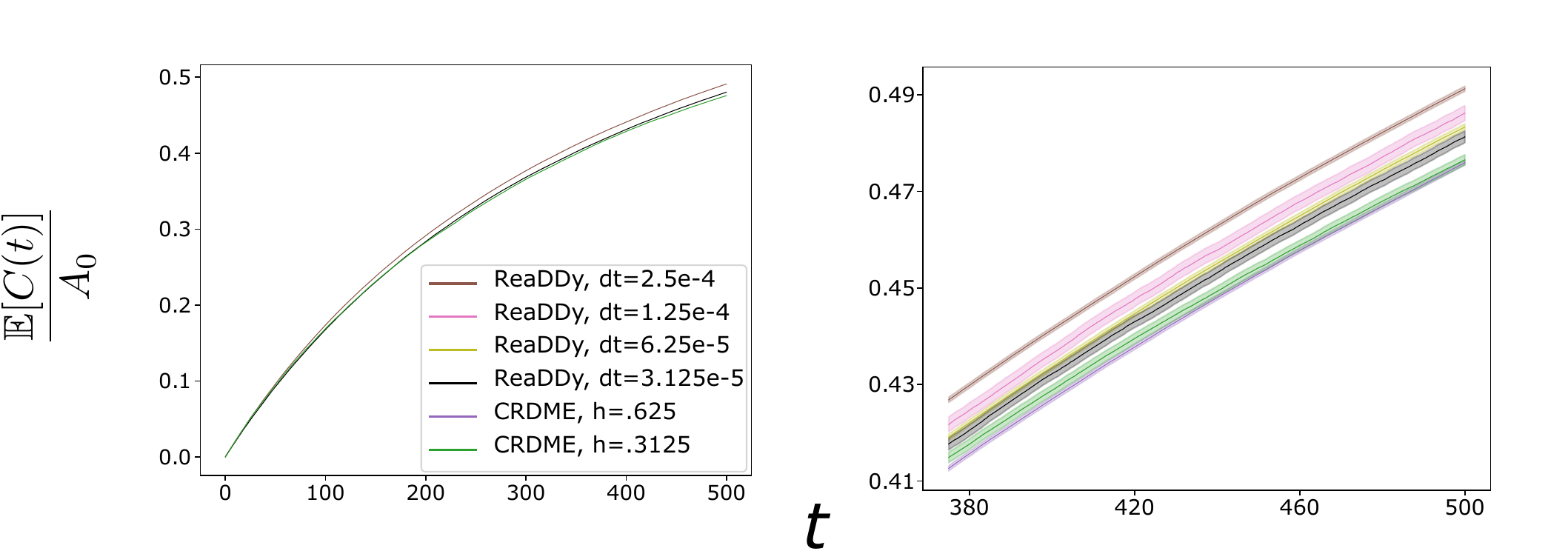}%
    \caption{Comparison of CRDME and ReaDDy transients for $A_0 = 800$ under different levels of mesh refinement and different time steps. The two CRDME simulation curves depicted only differ slightly, whereas the ReaDDy simulations steadily decrease toward the CRDME curves. In the left panel, we have omitted the pink, yellow, and purple curves for visual clarity. We have also included 95\% confidence intervals in the right panel. The CRDME curves are averaged over $2{,}800$ ($h = .625$) and $987$ ($h = .3125$) trials, and the ReaDDy curves over $2{,}791$, $1{,}000$, $1{,}979$, and $496$ trials for $\Delta t = 2.5\cdot 10^{-4}$, $1.25\cdot 10^{-4}$, $6.25\cdot 10^{-5}$, and $3.125\cdot 10^{-5}$, respectively.}
    \label{fig:gamma800_zoomed_comparision}
\end{figure}

We observe in Figure~\ref{fig:multiparticle_transient_comparision} that the ReaDDy transients and CRDME transients, averaged over $P = 2{,}800$ simulations for each $A_0$ ($2{,}791$ for the ReaDDy simulation with $A_0 = 800$), match very closely for $A_0 \leq 400$. For $A_0 = 800$, the first case for which we see a significant change in the shape of the curve, there is a much larger discrepancy, indicating that either the CRDME or ReaDDy simulations (or both) are not yet converged. We explore this case further in Figure~\ref{fig:gamma800_zoomed_comparision}, where we plot the mean proportion of bound particles $P_{\text{bound}}(t)$ for both the CRDME and ReaDDy simulations (with 95\% confidence intervals in the right panel) under further refinements of $h$ and $\Delta t$, respectively.

We observe in the left and right panels of Figure~\ref{fig:gamma800_zoomed_comparision} that the mean proportion of $C$ particles in the ReaDDy simulations steadily decreases toward that in the CRDME simulations as $\Delta t\to 0$. The trend is shown most clearly in the right panel, where we can also see that the CRDME curves for $h=.625,.3125$ are close with overlapping or nearly overlapping confidence intervals.

On the finest levels of refinement for CRDME and ReaDDy, each simulation required to produce the curves in Figure~\ref{fig:gamma800_zoomed_comparision} took 24--48 computing hours to complete. Because of limited computing time, we were not able to refine $h$ and $\Delta t$ further. However, we believe based on the results in Figure~\ref{fig:gamma800_zoomed_comparision} that continued refinement of $\Delta t$ would cause the ReaDDy solution to continue converging toward the CRDME solution.

\begin{figure}[h]
    \centering
    \includegraphics[width=.8\textwidth]{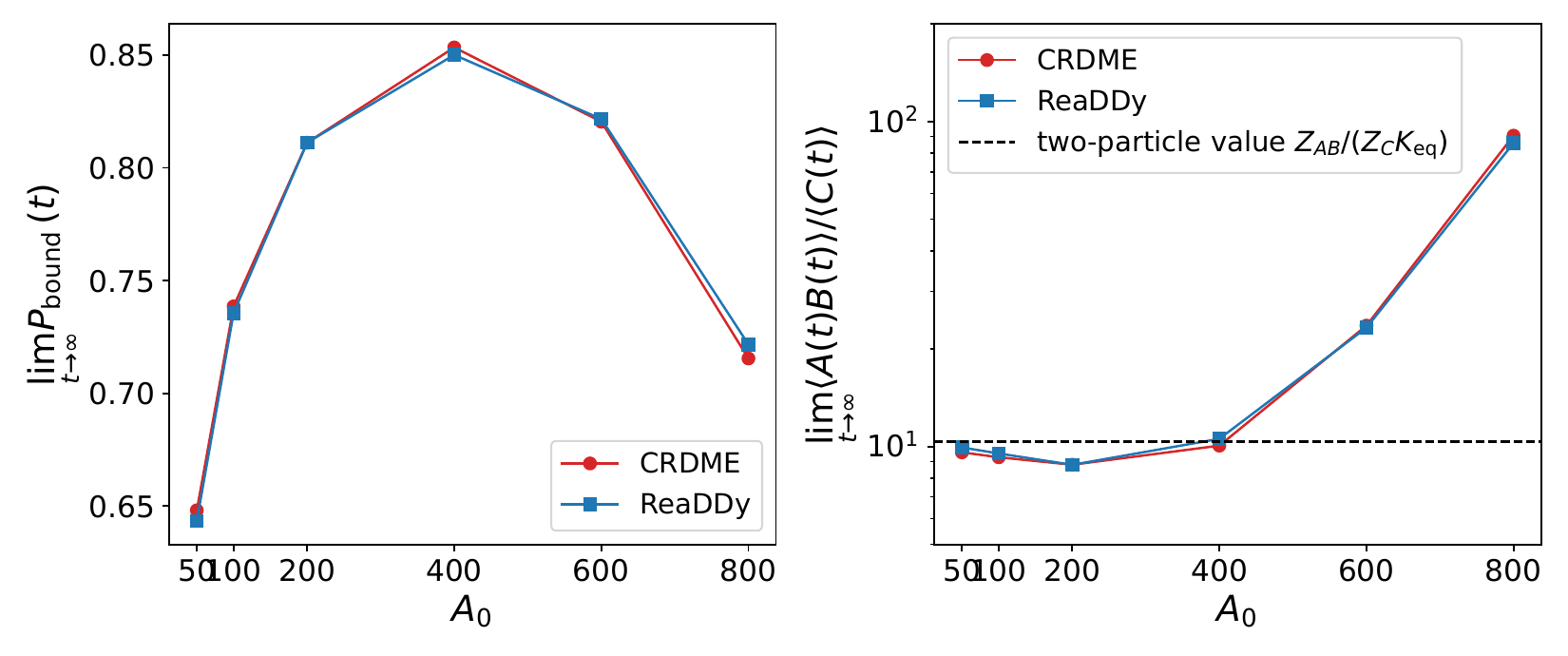}%
    \caption{Empirical values of the equilibrium bound fraction $\lim_{t\to\infty} P_{\text{bound}}(t)$ (left) and $\lim_{t\to\infty} \frac{\avg{A(t)B(t)}}{\avg{C(t)}}$ (right) for ReaDDy and CRDME, with particle counts $A_0=50,100,200,400,600,800$ (horizontal axis). The dashed black line on the right plot shows the analytical two-particle value.}
    \label{fig:KdV-multiparticle}
\end{figure}

Recall that $\bs_c = (M_0-c,M_0-c,c)$ is the number state with $c$ particles of type $C$, and note that $Z_C = \abs{\Omega}$ as there is no one-body potential for species $C$. In addition to the transients, we also compared the equilibrium ratios 
\begin{equation}\label{eq:equilibrium_statistic_multiparticle}
  \lim_{t\to\infty} \frac{\avg{A(t)B(t)}}{\avg{C(t)}}
  = \frac{\sum_{c = 0}^{M_0} (M_0-c)^2 \bar{\pi}_{\bs_c}}{\sum_{c=0}^{M_0} c\bar{\pi}_{\bs_c}}
  = \frac{1}{\Kd}\frac{\sum_{c=1}^{M_0} c\frac{Z_{\bs_{c-1}}}{Z_{\bs_{c}}} \bar{\pi}_{\bs_c}}{\sum_{c=1}^{M_0} c\bar{\pi}_{\bs_c}},
\end{equation} 
as we varied the initial condition $A_0 = M_0$. We derive the representation on the right-hand side of \eqref{eq:equilibrium_statistic_multiparticle} by using the relation \eqref{eq:multiparticle_DB_relation_partition_function_pis}. The difference between the plotted statistic and the two-particle version 
\begin{equation*}
  \lim_{t\to\infty} \frac{\avg{A(t)B(t)}}{\avg{C(t)}} = \frac{\bar{\pi}_{AB}}{\bar{\pi}_C} = \frac{Z_{AB}}{{Z_C}\Kd} \approx 10.3706
\end{equation*} 
(dashed black line in the right panel of Figure~\ref{fig:KdV-multiparticle}) reflects the divergence of $\frac{Z_{\bs_{c-1}}}{Z_{\bs_c}}$ from $\frac{Z_{AB}}{Z_C}$ at equilibrium as the number of particles increases and the system becomes more crowded; since the product $C$ has a larger radius than either substrate ($r_C = 1.4$ versus $r_A = r_B = 1$), crowding penalizes the bound state relative to unbound pairs, consistent with the steep increase of the plotted ratio for $A_0 \geq 600$; for $A_0 \leq 400$ the ratio stays close to the two-particle value, dipping slightly below it near $A_0 = 200$.

The steady-state values in Figure~\ref{fig:KdV-multiparticle} were obtained, for each code and each $A_0$, from ten independent trajectories run to time $t = 10{,}000$ and time-averaged after $t = 2{,}000$ ($t = 4{,}000$ for $A_0 = 800$). The results in Figure~\ref{fig:KdV-multiparticle} show good agreement between the CRDME simulations ($h = .625$) and the ReaDDy simulations ($\Delta t = 2.5\cdot 10^{-4}$). We note that we have used the same values of $h$ and $\Delta t$ that produced diverging transients in Figures~\ref{fig:multiparticle_transient_comparision} and \ref{fig:gamma800_zoomed_comparision} for $A_0 = 800$, showing that the curves eventually reconverge at steady state.

%% file: manuscript/derivation-general.tex
\section{Hopping and reaction rate derivations}\label{section:derivations}
 
 In this section, we will derive a discretization of \eqref{eq:FKE-full} for the model $A + B \rightleftharpoons C$ reaction from Section~\ref{section:introduction_reversible}, obtaining the CRDME approximation and the associated spatial hopping and reactive transition rates presented in Section~\ref{section:simulation_algorithm_pair_potentials} (summarized in Table~\ref{tab:discrete_reaction_rates}).
 
 We assume that we are given a triangulation $\mathcal{T}_h$ with vertices $\Omega_h = \{x_i\}_{i=1}^N$ and dual mesh voxels $\{V_i\}_{i=1}^N$ satisfying \eqref{eq:voxel-identities}, as before. 
 
 Denote by $V_h = \Span\{\phi_i\}_{i=1}^N$ the space generated by standard $P^1$ finite element basis functions $\phi_i$ associated with the nodes of $\mathcal{T}_h$. As before, $\bs = (s_1,\dots,s_M)$ will indicate the number state vector, with $s_m$ denoting the number of particles of species $S_m$, and $\mathbf{i}^{\bs} = (\mathbf{i}^{s_1},\dots,\mathbf{i}^{s_M}) \in \{1,\dots,N\}^{|\bs|}$ the vector of nodal positions for the particles (as in Section~\ref{section:simulation_algorithm_pair_potentials_hopping_rates}). Finally, we define the following product-space analogues of the single-particle mesh notations, associated with the multiparticle state vector $\mathbf{i}^{\bs}$:
 \begin{align*}
     \mathcal{T}_h^{\bs} &=  \left\{\bigtimes_{m=1}^M \bigtimes_{r=1}^{s_m} T^m_r : T^m_r \in \mathcal{T}_h \right\}, 
     & V_{\mathbf{i}^{{\bs}}} &= \bigtimes_{m=1}^M \bigtimes_{r=1}^{s_m} V_{i^{s_m}_r},\\
     \Omega_h^{\bs} &= \{\mathbf{x}_{\mathbf{i}^{\bs}}\},
     &\Phi^{\mathbf{i}^{\bs}}(\vq^{\bs}) &= \prod_{m=1}^M \prod_{r=1}^{s_m}\phi_{i^{s_m}_r}(q^{s_m}_r),&V_h^{\bs} &=  \Span\{\Phi^{\mathbf{i}^{\bs}}\}.
 \end{align*}

\subsection{Hopping rate derivation}\label{section:derivations_hopping}

Here we will derive the spatial hopping transition rates \eqref{eq:discrete_hopping_rate_formula} by discretizing the transport operator from \eqref{eq:FKE-full} in the case of a general multiparticle potential $\Psi^{\bs} : \Omega^{\abs{\bs}} \to \R$. For notational simplicity, we derive the discretization for a system with two species, $A$ and $B$, but note the same argument can be applied to a general multi-species system. In this special case we will take $\bs = (a,b)$ to denote the number state of $A$ and $B$ species, and $\bq^{\bs} = (\bq^{a},\bq^b)$ to denote the particle positions when in number state $\bs$. Ignoring the reaction terms, the Galerkin approximation of the weak form of \eqref{eq:FKE-full}, posed over the subspace $V_h^{\bs}\subseteq H^1(\Omega^{|\bs|})$, is to find $\density{\bs}_h\in V_h^{\bs}$ such that 
\begin{align}\label{eq:discrete-FKE-AB}
    \int_{\Omega^{|\bs|}} \frac{\partial \density{\bs}_h}{\partial t}(\bq^{\bs},t)v_h^{\bs} (\bq^{\bs})d\bq^{\bs} = a^{\bs}(\density{\bs}_h,v_h^{\bs}), 
    \quad \text{for all } v_h^{\bs} \in V_h^{\bs}, \, t \in[0,T],
\end{align}
where 
\begin{equation*}
   a^{\bs}(\density{\bs}_h,v_h^{\bs}) = -\int_{\Omega^{|\bs|}} 
   \begin{pmatrix} D^AI^{a d}& \\ 
        & D^BI^{b d} 
    \end{pmatrix}
    \left[\nabla \density{\bs}_h(\bq^{\bs},t) + \density{\bs}_h(\bq^{\bs},t)\nabla \Psi^\bs(\bq^{\bs})\right] \cdot \nabla v_h^{\bs}(\bq^{\bs}) \, d\bq^{\bs}.
\end{equation*}

We recall here the mass-lumping quadrature rule stemming from \eqref{eq:voxel-identities}, where we integrate $g\in C(\overline{\Omega})$ by replacing it with its nodal interpolant $g_h$:
\begin{equation} \label{eq:mass-lump-quadrature}
\int_\Omega g(x)dx \approx \mathcal{Q}^h_{\text{lump}}(g) := \int_\Omega g_h(x) dx = \sum_{i=1}^N g(x_i)|V_i|.
\end{equation}
One can show that for $g$ sufficiently regular, the mass-lumping quadrature has an error of order $h^2$. This follows since with $\mathcal{Q}^h_{\text{lump}}(g) = \int_\Omega g_h(x)\,dx$, the standard interpolation estimate \cite[Theorem~3.1.6]{ciarlet_finite_2002} gives $\left|\int_\Omega (g - g_h)\,dx\right| \leq |\Omega|^{1/2}\|g - g_h\|_{L^2(\Omega)} \leq C h^2 |g|_{H^2(\Omega)}$ for $g \in H^2(\Omega)$ and $d \leq 3$. See also \cite{raviart_use_1975,thomee_galerkin_2006} for the use of this quadrature in the lumped mass method for parabolic problems.
Although we do not prove this here, our numerical convergence studies support our expectation that the additional \textit{consistency error} introduced by mass-lumping quadrature does not change the convergence rate of our finite element discretization in either the $H^1$ or $L^2$ norm.  See, for example, \cite{ciarlet_finite_2002} for techniques that could be used to try to prove this, and~\cite{Heldman_2024} for a convergence proof of a related PDE discretization.

In what follows we denote the nodal values of $v_h^{\bs}\in V^{\bs}_h$ by $v_h^{\bs}(\node{\multinode{i}{s}}) =: v_{\mathbf{i}^\bs}$, and similarly write $\density{\bs}_{\multinode{i}{s}}(t) := \density{\bs}_h(\node{\multinode{i}{s}},t)$. We apply the mass-lumping quadrature rule \eqref{eq:mass-lump-quadrature} iteratively to the left-hand side of \eqref{eq:discrete-FKE-AB} to obtain 
\begin{align} \label{eq:crdme_lhs}
      \int_{\Omega^{|\bs|}} \frac{\partial \density{\bs}_h}{\partial t}(\bq^{\bs},t)v_h^{\bs} (\bq^{\bs})\,d\bq^{\bs} \approx \sum_{\multinode{i}{s}} |V_{\multinode{i}{s}}| \frac{\partial \density{\bs}_h}{\partial t}(\mathbf{x}_{\multinode{i}{s}},t)v_h^{\bs} (\mathbf{x}_{\multinode{i}{s}}) = \sum_{\multinode{i}{s}} |V_{\multinode{i}{s}}| \frac{d \density{\bs}_{\multinode{i}{s}}}{d t}(t)v_{\multinode{i}{s}} = \sum_{\multinode{i}{s}} \frac{d\cdf{\bs}_{\multinode{i}{s}}}{dt}v_{\multinode{i}{s}},
\end{align}
where we define $F_{\multinode{i}{s}}^{\bs}(t) := |V_{\multinode{i}{s}}|\density{\bolds}_{\multinode{i}{s}}(t)$. Note that $F^{\bs}_{\multinode{i}{s}}(t)$ is normalized so that $\tfrac{1}{\bs!}\sum_{\multinode{i}{s}} F^{\bs}_{\multinode{i}{s}}(t)$ is the probability of being in number state $\bs$ for all $t \in [0,T]$ and hence 
\begin{equation} \label{eq:crdme_normalization}
    \sum_{\bs} \frac{1}{\bs!}\sum_{\multinode{i}{s}} F^{\bs}_{\multinode{i}{s}}(t) = 1.
\end{equation}
While technically $F^{\bs}_{\multinode{i}{s}}(t) / \bs!$ is the probability of being in the state $(\bs,\multinode{i}{s})$, we will subsequently also refer to $F^{\bs}_{\multinode{i}{s}}(t)$ as a probability for ease of exposition.

We now consider an approximation of the right-hand side of \eqref{eq:discrete-FKE-AB}, the discrete spatial transport operator $a^{\bs}(\density{\bs}_h, v_h^{\bs})$. For the following discussion, we will drop the superscripts $\bs$, since spatial transport does not change the particle counts. We start by writing the full expression for $a(\density{}_h, v_h)$:
\begin{equation}\label{eq:FKE-AB-weak-form-rhs}
\begin{split}
    a(\density{}_h, v_h) = &-\sum_{\alpha=1}^a \int_{\Omega^{|\bs|}} D^A \left[\nabla_{q^a_\alpha}\density{}_h(\bq,t) + \density{}_h(\bq,t)\nabla_{q^a_\alpha}\Psi(\bq)\right]\cdot \nabla_{q^a_\alpha} v_h(\bq)d\bq
    \\&- \sum_{\beta=1}^b \int_{\Omega^{|\bs|}} D^B \left[\nabla_{q^b_\beta}\density{}_h(\bq,t) + \density{}_h(\bq,t)\nabla_{q^b_\beta}\Psi(\bq)\right]\cdot \nabla_{q^b_\beta} v_h(\bq) d\bq.
\end{split}
\end{equation}
We see how \eqref{eq:FKE-AB-weak-form-rhs} can be viewed as a sum of transport operators for each particle, making it tempting to discretize each term in the sum individually, as typically done in the absence of two-body interactions, see Section~\ref{section:simulation_algorithm_background_drift_diffusion}. Because $\Psi$ couples the coordinates, however, each term still depends on the positions of all particles, and the key idea of our discretization is to treat the coordinates within each term differently. We apply the mass-lumping quadrature rule with respect to all variables except the one corresponding to the particle whose motion is being described (i.e., the variable appearing in the gradients of a given term), and handle the latter using the EAFE method with linear nodal interpolation of the potential, as also described in Section~\ref{section:simulation_algorithm_background_drift_diffusion}.

It suffices to consider one of the terms in the sum from \eqref{eq:FKE-AB-weak-form-rhs}. To that end, let $\alpha \in \{1,2,\dots,a\}$ be one of the $A$ particles, and set 
\begin{align*}
\multinode{i}{s} \setminus i^a_\alpha &= (i^a_1,\dots,i^a_{\alpha-1},i^a_{\alpha+1},\dots,i^a_a,\mathbf{i}^b),\\   
(\multinode{i}{s} \setminus i^a_\alpha) \cup i &= (i^a_1,\dots,i^a_{\alpha-1},i,i^a_{\alpha+1},\dots,i^a_a,\mathbf{i}^b),\\
u_{\multinode{i}{s} \setminus i^a_\alpha}(x) &= u(x_{i^a_1},\dots,x_{i^a_{\alpha-1}},x,x_{i^a_{\alpha+1}},\dots,x_{i^a_a},\mathbf{x}_{\mathbf{i}^b}), 
\end{align*}
for any multi-index $\multinode{i}{s}$ and any $u:\Omega^{|\bs|}\to \RR$. For brevity, when setting $\bk = \multinode{i}{s} \setminus i^a_\alpha$ we will also use the notation $u_{\bk}(x) := u_{\multinode{i}{s} \setminus i^a_\alpha}(x)$ and $u_{\bk \cup j} := u_{(\multinode{i}{s} \setminus i^a_\alpha) \cup j} := u_{\bk}(x_j)$. Finally, for a function $u : \Omega \to \R$, we will use the difference operator notation
\begin{equation*}
    \bd_{ij} u := u(x_j) - u(x_i).
\end{equation*} 

Consider the $\alpha$th term in the sum over $A$ particles in the first line of~\eqref{eq:FKE-AB-weak-form-rhs}. Below, sums over $\bk := \bj^{\bs} \setminus j_\alpha^a$ run over all possible positions of the particles other than the $\alpha$th $A$ particle. Applying first the mass-lumping quadrature formula in every coordinate except $q^a_\alpha$, and then applying the EAFE quadrature formula \eqref{eq:EAFE_bilinear_form} to the remaining one-particle integrals, yields 
\begin{equation}
\label{eq:FKE-AB-binlinear-form}
\begin{aligned}
\int_{\Omega^{|\bs|}} &D^A \left[\nabla_{q^a_\alpha}\density{}_h(\bq,t) + \density{}_h(\bq,t)\nabla_{q^a_\alpha} \Psi(\bq)\right]\cdot \nabla_{q^a_\alpha} v_h(\bq)d\bq\\
&\approx \sum_{\bk := \multinode{j}{s} \setminus j^a_\alpha} |V_{\bk}|\int_\Omega D^A \left[ \nabla_{q^a_\alpha} \density{}_{\bk}(q^a_\alpha,t) + \density{}_{\bk}(q^a_\alpha,t)\nabla_{q^a_\alpha}\Psi_{\bk}(q^a_\alpha)\right]\cdot \nabla_{q^a_\alpha} v_{\bk}(q^a_\alpha)dq^a_\alpha\\
&\approx \sum_{\bk := \multinode{j}{s} \setminus j^a_\alpha} |V_{\bk}| \sum_{E_{ij}\in\mathcal{T}_h} D^A \omega_{ij}B\left(\bd_{i j} \Psi_{\bk}\right) \left(\density{}_{\bk \cup j}(t)e^{\bd_{i j}\Psi_{\bk}} - \density{}_{\bk\cup i}(t)\right)\bd_{i j}v_{\bk}. 
\end{aligned}
\end{equation}
Now, suppose that we choose $v_h = \Phi^{\multinode{i}{s}}$. Then, since 
\begin{equation*}
    \Phi^{\multinode{i}{s}}_{\multinode{j}{s} \setminus j^a_\alpha}(x_i) 
    = \delta_{\multinode{i}{s} \setminus i^a_\alpha, \multinode{j}{s} \setminus j^a_\alpha}\delta_{ii^a_\alpha} 
    = \delta_{\multinode{i}{s} \setminus i^a_\alpha, \multinode{j}{s} \setminus j^a_\alpha} \phi_{i^a_\alpha}(x_i)
\end{equation*}
for any $x_i \in \Omega_h$, proceeding from \eqref{eq:FKE-AB-binlinear-form} we have 
\begin{equation}
    \label{eq:FKE-AB-bilinear-form-single-row-variable-change}
\begin{aligned}
    \sum_{\bk := \multinode{j}{s} \setminus j^a_\alpha} &|V_{\bk}| \sum_{E_{ij}\in\mathcal{T}_h}  D^A\omega_{ij}B\left(\bd_{ij} \Psi_{\bk}\right) \left(\density{}_{\bk \cup j}(t)e^{\bd_{ij}\Psi_{\bk}} - \density{}_{\bk\cup i}(t)\right)\bd_{ij}
    \Phi^{\multinode{i}{s}}_{\bk} \\
    &= |V_{\multinode{i}{s} \setminus i^a_\alpha}| \sum_{E_{ij}\in\mathcal{T}_h} D^A\omega_{ij}B\left(\bd_{ij}\Psi_{\multinode{i}{s} \setminus i^a_\alpha}\right) \left(\density{}_{(\multinode{i}{s} \setminus i^a_\alpha)\cup j}(t)e^{\bd_{ij}\Psi_{\multinode{i}{s} \setminus i^a_\alpha}} - \density{}_{(\multinode{i}{s} \setminus i^a_\alpha)\cup i}(t)\right)\bd_{ij}\phi_{i^a_\alpha}\\
    &= |V_{\multinode{i}{s} \setminus i^a_\alpha}| D^A \bigg[ \sum_{E_{i i_{\alpha}^a}\in \mathcal{T}_h} \omega_{i i^a_\alpha}B\left(\bd_{{i i^{a}_\alpha}}\Psi_{\multinode{i}{s} \setminus i^a_\alpha}\right) \left(\density{}_{\multinode{i}{s}}(t)e^{\bd_{{i i^a_\alpha}}\Psi_{\multinode{i}{s} \setminus i^a_\alpha}} - \density{}_{(\multinode{i}{s} \setminus i^a_\alpha)\cup i}(t)\right) \\
    &\phantom{=} \qquad\qquad\qquad - \sum_{E_{i^a_\alpha j} \in \mathcal{T}_h}\omega_{i^a_\alpha j}B\left(\bd_{{i^{a}_\alpha j}}\Psi_{\multinode{i}{s} \setminus i^a_\alpha}\right) \left(\density{}_{(\multinode{i}{s} \setminus i^a_\alpha)\cup j}(t)e^{\bd_{i^a_\alpha j}\Psi_{\multinode{i}{s} \setminus i^a_\alpha}} - \density{}_{\multinode{i}{s}}(t)\right)\bigg]\\
    &= |V_{\multinode{i}{s} \setminus i^a_\alpha}| D^A \sum_{i \neq i_{\alpha}^a} \omega_{i i^a_\alpha} \left(B\left(\bd_{{i^{a}_\alpha} i}\Psi_{\multinode{i}{s} \setminus i^a_\alpha}\right) \density{}_{\multinode{i}{s}}(t)- B\left(\bd_{{i i^{a}_\alpha}}\Psi_{\multinode{i}{s} \setminus i^a_\alpha}\right) \density{}_{(\multinode{i}{s} \setminus i^a_\alpha)\cup i}(t)\right) \\
    &= \sum_{i \neq i_{\alpha}^a}  \left(\frac{D^A \omega_{i i^a_\alpha}}{|V_{i_{\alpha}^a}|} B\left(\bd_{{i^{a}_\alpha} i}\Psi_{\multinode{i}{s} \setminus i^a_\alpha}\right) \cdf{}_{\multinode{i}{s}}(t)- \frac{D^A \omega_{i^a_\alpha i}}{|V_i|} B\left(\bd_{{i i^{a}_\alpha}}\Psi_{\multinode{i}{s} \setminus i^a_\alpha}\right) \cdf{}_{(\multinode{i}{s} \setminus i^a_\alpha)\cup i}(t)\right),    
\end{aligned}
\end{equation}
where we have used that $B(-z) = e^z B(z)$ and that the stiffness matrix is symmetric, i.e., $\omega_{i j} = \omega_{j i}$, in obtaining the second-to-last equality.

Repeating the computation in~\eqref{eq:FKE-AB-bilinear-form-single-row-variable-change} for each $A$ and $B$ particle and substituting $\Phi^{\multinode{i}{s}}$ into~\eqref{eq:crdme_lhs} for $v_h^{\bs}$, we find that \eqref{eq:discrete-FKE-AB} reduces to the final master equation for the CTMC approximation $\cdf{}_{\multinode{i}{s}}(t)$. It is given by the coupled system of ODEs
\begin{equation}\label{eq:FKE-AB-generator}
\begin{split}
    \frac{d\cdf{}_{\multinode{i}{s}}}{dt} &= \sum_{\alpha = 1}^a
    \sum_{i \neq i_{\alpha}^a}  \left(\frac{D^A \omega_{i^a_\alpha i}}{|V_i|} B\left(\bd_{{i i^{a}_\alpha}}\Psi_{\multinode{i}{s} \setminus i^a_\alpha}\right) \cdf{}_{(\multinode{i}{s} \setminus i^a_\alpha)\cup i}(t) - \frac{D^A \omega_{i i^a_\alpha}}{|V_{i_{\alpha}^a}|} B\left(\bd_{{i^{a}_\alpha} i}\Psi_{\multinode{i}{s} \setminus i^a_\alpha}\right) \cdf{}_{\multinode{i}{s}}(t)\right) \\
    &\phantom{=} + \sum_{\beta = 1}^b \sum_{i \neq i_{\beta}^b}  \left(\frac{D^B \omega_{i^b_\beta i}}{|V_i|} B\left(\bd_{{i i^{b}_\beta}}\Psi_{\multinode{i}{s} \setminus i^b_\beta}\right) \cdf{}_{(\multinode{i}{s} \setminus i^b_\beta)\cup i}(t) - \frac{D^B \omega_{i i^b_\beta}}{|V_{i_{\beta}^b}|} B\left(\bd_{{i^{b}_\beta} i}\Psi_{\multinode{i}{s} \setminus i^b_\beta}\right) \cdf{}_{\multinode{i}{s}}(t)\right), 
\end{split}
\end{equation}
over all possible valid $\multinode{i}{s}$ and particle number states, $\bs$. We can identify the spatial hopping (i.e., transition) rates \eqref{eq:discrete_hopping_rate_formula} for the CTMC as the coefficients of this master equation. 

The master equation \eqref{eq:FKE-AB-generator} determines an associated spatial transition rate matrix $L_h^{\bs}$, and can be rewritten in the form \eqref{eq:discrete-FKE-full} with $R^{\bs}_h = 0$.  To define a valid CTMC, we need that $L_h^{\bs}$ has nonnegative off-diagonal entries and zero column sums. The former is satisfied whenever $\omega_{ij} \geq 0$ for all $i\neq j$, for example on Delaunay meshes in two dimensions, see \cite{xu_monotone_1999}, just as in the purely diffusive case~\cite{engblom_simulation_2009}. $L_h^{\bs}$ has zero column sums if and only if the constant vector is in its left nullspace. For a general finite element stiffness matrix $B_h$ with entries $b_h(\phi_j,\phi_i)$, one sees that the constant vector is in the left nullspace if and only if 
\begin{equation*}
    b_h(\phi_j, 1) = \sum_{i=1}^N b_h(\phi_j, \phi_i) = \sum_{i=1}^N [B_h]_{ij} = 0,\qquad j=1,2,\dots,N.
\end{equation*}
If $v_h$ is a constant function in \eqref{eq:FKE-AB-binlinear-form} then $\bd_{{ij}}v_{\multinode{i}{s} \setminus i^a_\alpha} = 0$ for all $i,j$. Hence the corresponding stiffness matrix has zero column sums.  Since $L_h^{\bs}$ and the stiffness matrix are connected by right-multiplication of the latter by the inverse diagonal lumped mass matrix, i.e., by rescaling the columns with the reciprocal voxel measures $|V_{\multinode{i}{s}}|^{-1}$, it follows that $L_h^{\bs}$ has zero column sums.

Finally, we now show that the generated master equation / CTMC model is consistent with detailed balance of the CTMC's spatial transitions holding at equilibrium (i.e., reversibility of the CTMC), and that it supports a Gibbs-Boltzmann equilibrium distribution. Let $\denseq_h\in V_h^{\bs}$ be the (unnormalized) equilibrium discrete Gibbs-Boltzmann distribution 
\begin{equation*}
    \denseq_h(\mathbf{x}_{\mathbf{i}^\bs}) = e^{-\Psi^{\bs}(\mathbf{x}_{\bi^\bs})},
\end{equation*}
and $\bar{F}^{{\bs}}_{\multinode{i}{s}} = \denseq_h(\mathbf{x}_{\mathbf{i}^\bs})|V_{\multinode{i}{s}}|.$ Consider two states that differ only in the position of one $A$ particle, $\multinode{i}{s}$ and $\multinode{j}{s},$ where $\multinode{j}{s} = (\multinode{i}{s} \setminus i^a_\alpha) \cup j$ with vertex $j \in \{1,\dots,N\}$ and $j \neq {i^a_\alpha}$. Note the identities $\denseq_h(\mathbf{x}_{\multinode{i}{s}}) = \denseq_h(\mathbf{x}_{\multinode{j}{s}})\,e^{\bd_{i^a_\alpha j} \Psi_{\bi^\bs \setminus i_\alpha^a}}$, $e^{z}B(z) = B(-z)$ for all $z \in \RR$, $\frac{|V_{\multinode{i}{s}}|}{|V_{i^a_\alpha}|} = \frac{|V_{\multinode{j}{s}}|}{|V_{j}|}$, and $\omega_{i j} = \omega_{j i}$. Then
\begin{align*}
[L^{\bs}_h]_{\multinode{j}{s}\multinode{i}{s}} \cdfeq_{\multinode{i}{s}}  &= 
\frac{D^A \omega_{j i^a_\alpha}}{|V_{i_{\alpha}^a}|} B\left(\bd_{i^{a}_\alpha j}\Psi_{\multinode{i}{s} \setminus i^a_\alpha}\right) \cdfeq_{\multinode{i}{s}} 
= \frac{D^A \omega_{i^a_\alpha j}|V_{\bi^\bs}|}{|V_{i_{\alpha}^a}|} B\left(-\bd_{j {i^{a}_\alpha}}\Psi_{\multinode{i}{s} \setminus i^a_\alpha}\right) \denseq_h(\mathbf{x}_{\multinode{i}{s}}) \\
&= \frac{D^A \omega_{i^a_\alpha j}|V_{\bj^\bs}|}{|V_{j}|} B\left(\bd_{j {i^{a}_\alpha}}\Psi_{\multinode{i}{s} \setminus i^a_\alpha}\right) e^{-\bd_{{i^{a}_\alpha} j}\Psi_{\multinode{i}{s} \setminus i^a_\alpha}} \denseq_h(\mathbf{x}_{\multinode{i}{s}})
= \frac{D^A \omega_{i^a_\alpha j}}{|V_{j}|} B\left(\bd_{j {i^{a}_\alpha}}\Psi_{\multinode{i}{s} \setminus i^a_\alpha}\right) \cdfeq_{\multinode{j}{s}}  = [L^{\bs}_h]_{\multinode{i}{s}\multinode{j}{s}} \cdfeq_{\multinode{j}{s}}.
\end{align*}
A similar argument can be made for the transition rates at equilibrium of $B$ particles to neighboring voxels, i.e., transitions from state $\mathbf{i}^{\bs}$ to $(\mathbf{i}^\bs \setminus i^b_\beta)\cup j$ for $j$ adjacent to $i^b_\beta$. It follows, then, that the discretized CRDME is consistent with a Gibbs-Boltzmann equilibrium distribution in the absence of reactions, and that the associated CTMC is reversible and satisfies detailed balance with respect to this distribution. 

Assume the mesh is connected and sufficiently regular that the graph formed by edges with $\omega_{ij} > 0$ is connected. Likewise, assume all potential differences are finite. The spatial transition rate matrix $L^{\bs}_h$ is then the generator of a finite, irreducible CTMC. $L^{\bs}_h$ then has a one-dimensional nullspace spanned by the discrete Gibbs-Boltzmann distribution $\bar{F}^{\bs}_{\multinode{i}{s}}$. In the special case that there are no reactions and the number state is exactly $\bs$,  $\bs! \bar{F}^{\bs}_{\multinode{i}{s}} / (\sum_{\bj^\bs} \bar{F}^{\bs}_{\multinode{j}{s}})$ gives the unique stationary (equilibrium) distribution of the CTMC under the CRDME normalization~\eqref{eq:crdme_normalization}.

\subsection{Reaction rate derivations}\label{section:derivations_reactions}

We now derive the general multiparticle CRDME for the reversible $A + B \rightleftharpoons C$ reaction, with reactive transition rates \eqref{eq:reaction-rates-general-form}. The volume average transition rates for the purely diffusive CRDME given in~\eqref{eq:discrete_fwd_standard} and~\eqref{eq:discrete_bwd_standard} are used as the base from which we derive modified rates that account for potential interactions.
The other reactive transition rate approximations, including \eqref{eq:discrete_fwd_lump} / \eqref{eq:discrete_bwd_lump} and \eqref{eq:discrete_fwd_smooth} / \eqref{eq:discrete_bwd_smooth}, follow from making additional lumping / quadrature approximations to the resulting integrals for the volume averages as described in Section~\ref{section:simulation_algorithm_background_reactions}.

 \subsubsection{Derivation of \eqref{eq:reaction-rates-general-form}}\label{section:derivations_reactions_multiparticle}

Consider the reaction operator $R^{\bs}$ from the full particle model given in~\eqref{eq:reversible-operator}. Let $M = 3$, let $\bs = (a,b,c) \in(\NN_0)^M$ be the particle number state, and let $\multinode{i}{s} \in \{1,2,\dots,N\}^{|\bs|}$. We define a new notation that will be useful at intermediate steps of our derivation. The annihilation operator $\annihilation{{i_1^{s_{k_1}}i_2^{s_{k_2}}}\dots i_K^{{s_{k_K}}}}\bq$ will return the vector of particle positions $\bq$ with each particle $i_\ell^{s_{k_\ell}}$ of species $S_{k_\ell}$ removed. 
Finally, recall that $V_{i j} = V_i \times V_j$ and $V_{i j k} = V_i \times V_j \times V_k$ denote product voxels.

With these notations in mind, we discretize each of the four terms in~\eqref{eq:reversible-operator} separately. After multiplication by $\Phi^{\multinode{i}{s}}$ and integration over $\Omega^{|\bs|}$, the first term is approximated by
\begin{align*}
    &\int_{\Omega^{|\bs|}}\density{\bs}(\bq,t) \left[\sum_{\alpha=1}^a\sum_{\beta=1}^b \int_\Omega \kappa^+(z|q^a_\alpha,q^b_\beta)\nu^+(\fwdop{\alpha}{\beta}{z}\boldq , \bq)dz\right]\Phi^{\multinode{i}{s}}(\bq)\,d\bq
    \\&\approx \sum_{\alpha=1}^a\sum_{\beta=1}^b \frac{|V_{\multinode{i}{s}}|}{|V_{i^a_\alpha i^b_\beta}|}\int_{\Omega^{|\bs|}} \density{\bs}(\bq,t) \left[\int_\Omega \kappa^+(z|q^a_\alpha,q^b_\beta)\nu^+(\fwdop{\alpha}{\beta}{z}\boldq , \bq)dz\right]\phi_{i^a_\alpha}(q^a_\alpha) \phi_{i^b_\beta}(q^b_\beta) \delta\left(\annihilation{i^a_\alpha i^b_\beta}\left[\bq - \mathbf{x}_{\multinode{i}{s}}\right]\right)d\bq\\
    &\approx \sum_{\alpha=1}^a\sum_{\beta=1}^b \sum_{k=1}^N \cdf{\bs}_{\multinode{i}{s}}(t)\frac{\nu^+({\fwdop{\alpha}{\beta}{k}\bx_{\multinode{i}{s}},\bx_{\multinode{i}{s}}})}{|V_{i^a_\alpha i^b_\beta}|}\int_{V_{i^a_\alpha i^b_\beta k}} \kappa^+(z|x,y)dzdxdy \\
    &= \brac{\sum_{\alpha=1}^a\sum_{\beta=1}^b \sum_{k=1}^N \hat{\kappa}^+_{i^a_\alpha i^b_\beta k} (\bi^\bs)}\cdf{\bs}_{\multinode{i}{s}}(t).
\end{align*}
To go from the first line to the second line, we apply mass-lumping quadrature to the integrals over all coordinates except for $q^a_\alpha$ and $q^b_\beta$. By \eqref{eq:mass-lump-quadrature}, this involves multiplication by the volume factors in $|V_{\multinode{i}{s}}| / |V_{i^a_\alpha i^b_\beta}|$  and evaluation of the integrands with respect to non-reactant coordinates at the mesh nodes. The latter is represented here by the Dirac delta function $\delta\left(\annihilation{i^a_\alpha i^b_\beta}\left[\bq - \mathbf{x}_{\multinode{i}{s}}\right]\right)$. To go from the second to the third line, we first use the multidimensional Dirac delta to integrate out all non-reactant coordinates in $\bq$ (i.e., all but $q_\alpha^a$ and $q_\beta^b$). We then
apply the same hybrid mass-lumping / finite volume quadrature underlying~\eqref{eq:discrete_fwd_standard} to the remaining integral over $z$, $q_\alpha^a$, and $q_\beta^b$. All factors within the integrand except $\kappa^+$ are evaluated at the corresponding pre- and post-reaction nodal configurations, while $\kappa^+$ is retained within the integral over each product voxel $V_{i_\alpha^a i_\beta^b k}$. 
In the absence of potentials, this is equivalent to the hybrid finite element--finite volume approach we used for reaction-diffusion problems, where reactive terms were handled via finite volume approximation, see~\cite{isaacson_unstructured_2018}. The last line follows by the definition of $\hat{\kappa}^+_{i^a_\alpha {i^b_\beta} k}$ in~\eqref{eq:reaction-rates-general-form}.

The discretization of the other three terms in \eqref{eq:reversible-operator} is similar; we first apply mass-lumping quadrature with respect to the variables representing particles not involved in the reaction (line 1 to line 2), and then apply mass-lumping quadrature to all terms for the particles involved in the reaction except $\kappa^{\pm}$, obtaining a finite volume average for the latter. For the other outgoing flux from the number state $\bs$ we obtain
\begin{align*}
    &\int_{\Omega^{|\bs|}} \density{{\bs}}(\qvec,t) \left[\sum_{\xi = 1}^c \int_{\Omega^2}\kappa^-(x,y|q^c_\xi)\nu^-(\bwdop{x}{y}{\xi} \qvec,\qvec)dxdy\right] \Phi^{\multinode{i}{s}}(\qvec)\,d\qvec\\
    &\approx \sum_{\xi = 1}^c \frac{|V_{\multinode{i}{s}}|}{|V_{i^c_\xi}|}\int_{\Omega^{|\bs|}} \density{\bs}(\qvec,t)  \left[\int_{\Omega^2}\kappa^-(x,y|q^c_\xi)\nu^-(\bwdop{x}{y}{\xi} \qvec,\qvec)dxdy\right] \phi_{i^c_\xi}(q^c_\xi) \delta\left(\annihilation{i^c_\xi}\left[\bq - \mathbf{x}_{\multinode{i}{s}}\right]\right) d\qvec\\
    &\approx \sum_{\xi=1}^c \sum_{i=1}^N \sum_{j=1}^N \cdf{\bs}_{\multinode{i}{s}}(t)\frac{\nu^-\paren{{\bwdop{i}{j}{\xi}\bx_{\multinode{i}{s}},\bx_{\multinode{i}{s}}}}}{|V_{i^c_\xi}|}\int_{V_{i^c_\xi ij}}\kappa^-(x,y|z)dzdxdy\\
    &= \brac{\sum_{\xi=1}^c \sum_{i=1}^N \sum_{j=1}^N \hat{\kappa}^-_{i_\xi^c i j}(\multinode{i}{s})} \cdf{\bs}_{\multinode{i}{s}}(t).
\end{align*}
For the two incoming fluxes we let $\bolds^+ = (a-1,b-1,c+1)$ and $\bolds^- = (a+1,b+1,c-1)$ for compactness and have 
\begin{align*}
    &\int_{\Omega^{|\bs|}}\left[\sum_{\xi=1}^c\int_{\Omega^2} \kappa^+(q^c_\xi | x,y)\nu^+(\bq,\bwdop{x}{y}{\xi}\qvec)\density{\bs^-}(\bwdop{x}{y}{\xi}\qvec,t)dxdy\right]\Phi^{\multinode{i}{s}}(\qvec)d\qvec
    \\&\approx \sum_{\xi=1}^c \frac{|V_{\multinode{i}{s}}|}{|V_{i^c_\xi}|}\int_{\Omega^{|\bs|}}\left[\int_{\Omega^2} \kappa^+(q^c_\xi | x,y)\nu^+(\bq,\bwdop{x}{y}{\xi}\qvec)\density{\bs^-}(\bwdop{x}{y}{\xi}\qvec,t)dxdy\right]\phi_{i^c_\xi}(q^c_\xi)\delta\left(\annihilation{i^c_\xi}\left[\bq - \mathbf{x}_{\multinode{i}{s}}\right]\right) d\qvec
    \\&\approx \sum_{\xi=1}^c \sum_{i=1}^N \sum_{j=1}^N \brac{\frac{\nu^+\paren{\bx_{\multinode{i}{s}},\bx_{\bwdop{i}{j}{\xi}\multinode{i}{s}}}}{|V_{i j}|}\int_{V_{i j i^c_\xi }} \kappa^+(z | x,y)dxdydz}\cdf{\bs^-}_{\bwdop{i}{j}{\xi}\multinode{i}{s}}(t) \\
    &= \sum_{\xi=1}^c \sum_{i=1}^N \sum_{j=1}^N \hat{\kappa}^+_{i j i_\xi^c}(\bwdop{i}{j}{\xi}\multinode{i}{s}) \cdf{\bs^-}_{\bwdop{i}{j}{\xi}\multinode{i}{s}}(t),
\end{align*}
and finally 
\begin{align*}
& \int_{\Omega^{|\bs|}} \left[ \sum_{\alpha=1}^a\sum_{\beta=1}^b \int_\Omega \kappa^-(q^a_\alpha,q^b_\beta | z)\nu^-(\bq , \fwdop{\alpha}{\beta}{z}\qvec )\density{\bs^+}(\fwdop{\alpha}{\beta}{z}\qvec,t)dz\right] \Phi^{\multinode{i}{s}}(\qvec)\, d\qvec
\\&\approx \sum_{\alpha=1}^a\sum_{\beta=1}^b  \frac{|V_{\multinode{i}{s}}|}{|V_{i^a_\alpha i^b_\beta}|}\int_{\Omega^{|\bs|}}\!\left[ \int_\Omega \kappa^-(q^a_\alpha,q^b_\beta | z)\nu^-(\bq , \fwdop{\alpha}{\beta}{z}\qvec )\density{\bs^+}(\fwdop{\alpha}{\beta}{z}\qvec,t)dz\right]\!\phi_{i^a_\alpha}(q^a_\alpha)\phi_{i^b_\beta}(q^b_\beta)\delta\!\left(\annihilation{i^a_\alpha i^b_\beta}\left[\bq - \mathbf{x}_{\multinode{i}{s}}\right]\right)d\qvec
\\&\approx \sum_{\alpha=1}^a\sum_{\beta=1}^b\sum_{k=1}^N \brac{\frac{\nu^-\paren{\bx_{\multinode{i}{s}},\bx_{\fwdop{\alpha}{\beta}{k}\multinode{i}{s}}}}{|V_{k}|}\int_{V_{ k i^a_\alpha i^b_\beta}} \kappa^-(x,y|z)dzdxdy}\cdf{\bs^+}_{\fwdop{\alpha}{\beta}{k}\multinode{i}{s}}(t) \\
&= \sum_{\alpha=1}^a\sum_{\beta=1}^b\sum_{k=1}^N \hat{\kappa}^-_{k i_\alpha^a i_\beta^b}(\fwdop{\alpha}{\beta}{k}\multinode{i}{s}) \cdf{\bs^+}_{\fwdop{\alpha}{\beta}{k}\multinode{i}{s}}(t).
\end{align*}
Collectively these define $R_h^{\bs}$ in~\eqref{eq:discrete-FKE-full}. 

The complete discretized CRDME approximation to the underlying particle model~\eqref{eq:FKE-full} for the reversible $A + B \rightleftharpoons C$ reaction, showing both spatial hops that model the drift-diffusion of particles between voxels and reactive interactions between particles in nearby voxels, is then
\begin{align*}
    \frac{d\cdf{\bolds}_{\multinode{i}{s}}}{dt}(t) 
    &= \sum_{\alpha = 1}^a
    \sum_{i \neq i_{\alpha}^a}  \left(\frac{D^A \omega_{i^a_\alpha i}}{|V_i|} B\left(\bd_{{i i^{a}_\alpha}}\Psi_{\multinode{i}{s} \setminus i^a_\alpha}\right) \cdf{\bs}_{(\multinode{i}{s} \setminus i^a_\alpha)\cup i}(t) - \frac{D^A \omega_{i i^a_\alpha}}{|V_{i_{\alpha}^a}|} B\left(\bd_{{i^{a}_\alpha} i}\Psi_{\multinode{i}{s} \setminus i^a_\alpha}\right) \cdf{\bs}_{\multinode{i}{s}}(t)\right) \\
    &\phantom{=} + \sum_{\beta = 1}^b \sum_{i \neq i_{\beta}^b}  \left(\frac{D^B \omega_{i^b_\beta i}}{|V_i|} B\left(\bd_{{i i^{b}_\beta}}\Psi_{\multinode{i}{s} \setminus i^b_\beta}\right) \cdf{\bs}_{(\multinode{i}{s} \setminus i^b_\beta)\cup i}(t) - \frac{D^B \omega_{i i^b_\beta}}{|V_{i_{\beta}^b}|} B\left(\bd_{{i^{b}_\beta} i}\Psi_{\multinode{i}{s} \setminus i^b_\beta}\right) \cdf{\bs}_{\multinode{i}{s}}(t)\right) \\
    &\phantom{=} + \sum_{\xi = 1}^c \sum_{i \neq i_{\xi}^c}  \left(\frac{D^C \omega_{i^c_\xi i}}{|V_i|} B\left(\bd_{{i i^{c}_\xi}}\Psi_{\multinode{i}{s} \setminus i^c_\xi}\right) \cdf{\bs}_{(\multinode{i}{s} \setminus i^c_\xi)\cup i}(t) - \frac{D^C \omega_{i i^c_\xi}}{|V_{i_{\xi}^c}|} B\left(\bd_{{i^{c}_\xi} i}\Psi_{\multinode{i}{s} \setminus i^c_\xi}\right) \cdf{\bs}_{\multinode{i}{s}}(t)\right) \\
    &\phantom{=} - \brac{\sum_{\alpha=1}^a\sum_{\beta=1}^b \sum_{k=1}^N \hat{\kappa}^+_{i^a_\alpha i^b_\beta k} (\bi^\bs)}\cdf{\bs}_{\multinode{i}{s}}(t)
    - \brac{\sum_{\xi=1}^c \sum_{i=1}^N \sum_{j=1}^N \hat{\kappa}^-_{i_\xi^c i j}(\multinode{i}{s})} \cdf{\bs}_{\multinode{i}{s}}(t)\\
    &+ \sum_{\xi=1}^c \sum_{i=1}^N \sum_{j=1}^N \hat{\kappa}^+_{i j i_\xi^c}(\bwdop{i}{j}{\xi}\multinode{i}{s}) \cdf{\bs^-}_{\bwdop{i}{j}{\xi}\multinode{i}{s}}(t)  
    + \sum_{\alpha=1}^a\sum_{\beta=1}^b\sum_{k=1}^N \hat{\kappa}^-_{k i_\alpha^a i_\beta^b}(\fwdop{\alpha}{\beta}{k}\multinode{i}{s}) \cdf{\bs^+}_{\fwdop{\alpha}{\beta}{k}\multinode{i}{s}}(t).
\end{align*}
This shows the CTMC reaction structure with transition rates~\eqref{eq:reaction-rates-general-form}, and more generally implies the transitions given in Table~\ref{tab:discrete_reaction_rates} for spatial hops and the $A_i + B_j \rightleftharpoons C_{{k}}$ reactions. It also identifies both $L_h^{\bs}$ and $R_h^\bs$ in~\eqref{eq:discrete-FKE-full}.

Finally, we now demonstrate that under the assumptions of Section~\ref{section:simulation_algorithm_background} on the reaction kernels of the underlying PBSRD model, the preceding CRDME is consistent with the principle of detailed balance for pointwise reactive fluxes holding at equilibrium (when in a bounded domain with reflecting boundary conditions). As described in~\cite{zhang_detailed_2022}, we expect detailed balance to hold for any reversible mass action reaction in a closed system from equilibrium statistical mechanics / thermodynamics considerations. It corresponds to the statement that the forward and backward reactive fluxes should balance at equilibrium. Let $k$ denote the voxel index for a newly produced $C$ particle following an $A + B \to C$ reaction between the $\alpha$th $A$ particle and the $\beta$th $B$ particle from the number state $\bs$ with particle positions $\bi^\bs$. Mathematically, the principle of detailed balance for this transition states that
\begin{equation}\label{eq:multpart_crdme_db}
\hat{\kappa}^+_{i^a_\alpha i^b_\beta k} (\bi^\bs)\bar{F}^{\bs}_{\multinode{i}{s}} = \hat{\kappa}^-_{k i_\alpha^a i_\beta^b}(\fwdop{\alpha}{\beta}{k}\multinode{i}{s}) \bar{F}^{\bs^+}_{\fwdop{\alpha}{\beta}{k}\multinode{i}{s}}.
\end{equation}
This implies that $[R_h^{\bs} \bar{\cdfvec}_h]_{\bi^{\bs}} = 0$ for all $\bs$ and $\bi^{\bs}$, which, since we are at steady state, in turn implies that $[L_h^{\bs} \bar{F}^{\bs}]_{\bi^{\bs}} = 0$ for all $\bs$ and $\bi^{\bs}$ by~\eqref{eq:discrete-FKE-full}. Under the same mesh and potential assumptions as in Section~\ref{section:derivations_hopping}, since the transition rate matrix $L_h^{\bs}$ is irreducible and finite, it will have a one-dimensional nullspace spanned by the discrete Gibbs-Boltzmann distribution for each fixed $\bs$. $\bar{F}^{\bs}_{\bi^{\bs}}$ can therefore be written in terms of a Gibbs-Boltzmann distribution as
\begin{equation*}
\bar{F}^{\bs}_{\bi^{\bs}} = \frac{\bs! \bar{\pi}_h^{\bs}}{Z_h^{\bs}} e^{-\Psi^{\bs}(\bx_{\bi^\bs})} |V_{\bi^\bs}|,
\end{equation*}
where $\bar{\pi}_h^{\bs}$ denotes the equilibrium probability of being in number state $\bs$ and $Z_h^\bs$ is the partition function for the discrete Gibbs-Boltzmann distribution given by
\begin{equation*}
    Z^{\bs}_h = \sum_{\bi^{\bs}} e^{-\Psi^{\bs}(\bx_{\bi^\bs})} |V_{\bi^\bs}|.
\end{equation*}
Substituting $\bar{F}^{\bs}_{\bi^{\bs}}$ into the discrete detailed balance relation~\eqref{eq:multpart_crdme_db}, we find the reactive transition rates are related by
\begin{equation} \label{eq:multpart_crdme_db_step2}
\hat{\kappa}^-_{k i_\alpha^a i_\beta^b}(\fwdop{\alpha}{\beta}{k}\multinode{i}{s})
= \hat{\kappa}^+_{i^a_\alpha i^b_\beta k} (\bi^\bs) \frac{\bar{F}^{\bs}_{\multinode{i}{s}}}{\bar{F}^{\bs^+}_{\fwdop{\alpha}{\beta}{k}\multinode{i}{s}}}
= \hat{\kappa}^+_{i^a_\alpha i^b_\beta k} (\bi^\bs) \frac{a b}{c+1}\frac{Z_h^{\bs^+}}{Z_h^{\bs}}\frac{|V_{i_{\alpha}^a}||V_{i_{\beta}^b}|}{|V_{k}|} \frac{\bar{\pi}_h^{\bs}}{\bar{\pi}_h^{\bs^+}} e^{\brac{\Psi^{\bs^+}\paren{\bx_{\fwdop{\alpha}{\beta}{k}\multinode{i}{s}}}-\Psi^{\bs}(\bx_{\bi^{\bs}})}}.
\end{equation}
Denote the implicit spatial constant in \eqref{eq:DB-ratio} by $\Theta^{\bs}$, i.e.,
\begin{equation*}
    \frac{\nu^-\paren{\bx_{\multinode{i}{s}},\bx_{\fwdop{\alpha}{\beta}{k}\multinode{i}{s}}}}{\nu^+\paren{\bx_{\fwdop{\alpha}{\beta}{k}\multinode{i}{s}},\bx_{\multinode{i}{s}}}} = \Theta^\bs e^{\brac{\Psi^{\bs^+}\paren{\bx_{\fwdop{\alpha}{\beta}{k}\multinode{i}{s}}}-\Psi^{\bs}(\bx_{\bi^{\bs}})}}.
\end{equation*}
Substituting the preceding, \eqref{eq:discrete-DB-condition-no-potentials}, and \eqref{eq:reaction-rates-general-form} into~\eqref{eq:multpart_crdme_db_step2} and rewriting, we obtain that for the discrete CRDME model
\begin{equation*}
\Theta^{\bs} = \Kd^h \frac{a b}{c+1}\frac{Z_h^{\bs^+}}{Z_h^{\bs}} \frac{\bar{\pi}_h^{\bs}}{\bar{\pi}_h^{\bs^+}}.
\end{equation*}
If we choose $\Theta^{\bs} = 1$, as was done in the underlying continuum PBSRD model for the concrete choices of modulating factors that have been proposed in practice, we obtain
\begin{equation} \label{eq:eq_pi_relation}
    \frac{c+1}{a b} \frac{\bar{\pi}_h^{\bs^+}}{\bar{\pi}_h^{\bs}} =  \Kd^h \frac{Z_h^{\bs^+}}{Z_h^{\bs}},
\end{equation}
which is analogous to the detailed balance constraint~\eqref{eq:multiparticle_DB_relation_partition_function_pis} from the continuum PBSRD model. 

Assume that we have a finite number of particles of each type, $A$, $B$, and $C$, in the system, and note the conservation laws $a + c = M_A$ and $b + c = M_B$, where $M_A$ and $M_B$ are the total numbers of $A$ and $B$ constituents when including those bound together. The set of number states, $\bs$, that are reachable for specific values of $M_A$ and $M_B$ is finite, and represents a finite stoichiometric compatibility class for the system. Equation~\eqref{eq:eq_pi_relation} then provides a recursive relation from which we can solve for the equilibrium probabilities $\bar{\pi}_h^{\bs}$ for all number states $\bs$ in the stoichiometric compatibility class, up to a common constant that can then be determined by normalization. Given these equilibrium probabilities $\bar{\pi}_h^{\bs}$, we have then established that $\bar{F}^{\bs}_{\bi^\bs}$ satisfies the detailed balance relationship~\eqref{eq:multpart_crdme_db}. Since detailed balance implies that $[R_h^{\bs}\bar{\cdfvec}_h]_{\bi^\bs} = 0$,   $\bar{F}^{\bs}_{\bi^\bs}$ was chosen to satisfy $[L_h^{\bs}\bar{F}^{\bs}]_{\bi^\bs} = 0$, and $\bar{F}^{\bs}_{\bi^\bs}$ is normalized per~\eqref{eq:crdme_normalization}, we conclude that $\bar{F}^{\bs}_{\bi^\bs}$ is an equilibrium solution to the CRDME within each stoichiometric compatibility class.

Under the mesh and potential assumptions of Section~\ref{section:derivations_hopping}, the mesh forms a connected graph (with respect to the edges with $\omega_{ij} > 0$), and the spatial transitions connect any two spatial configurations of particles for a fixed number state. Likewise, the reversible reactions connect any two consecutive number states in the chain representing the associated stoichiometric compatibility class.  The CTMC associated with the CRDME is therefore irreducible within each such class. Since it is also finite, there is a one-dimensional nullspace for the full transition rate matrix for the coupled CRDME over all $\bs$ in a given compatibility class, and the (normalized) equilibrium distribution $\bar{F}^{\bs}_{\bi^\bs}$ is unique within each stoichiometric compatibility class. 

For the full time-dependent CRDME with a fixed initial state $A(0) = a_0$, $B(0) = b_0$, and $C(0) = c_0$, which lies in a single stoichiometric compatibility class, the finiteness of that class and the irreducibility of the CTMC within it imply that the system converges as $t \to \infty$ to the unique equilibrium distribution of the class, given by $\bar{F}^{\bs}_{\bi^\bs}$. More generally, for an initial distribution spanning multiple compatibility classes, the probability assigned to each class is conserved, and the solution converges to the unique limiting equilibrium determined by these initial class probabilities, i.e., the corresponding mixture of the class-specific equilibrium distributions. 

In summary, we have established that the CRDME satisfies the principle of detailed balance at equilibrium. On each stoichiometric compatibility class, it defines a reversible CTMC with a conditional Gibbs-Boltzmann spatial equilibrium distribution. For each initial distribution, the associated time-dependent CRDME solution converges to a uniquely determined mixture of such class-specific distributions as $t \to \infty$.